\documentclass[11pt,twoside]{article}
\usepackage{latexsym}
\usepackage{amssymb,amsbsy,amsmath,amsfonts,amssymb,amscd}
\usepackage{graphicx,dsfont}
\usepackage{hyperref}
\usepackage{caption,color}
\newcommand{\commentout}[1]{}
\newcommand{\R}{\mathbb{R}}

\newcommand {\fer}   {\eqref}

\newcommand {\al} {\alpha}
\newcommand {\e}  {\varepsilon}
\newcommand {\da} {\delta}

\newcommand {\vp} {\varphi}

\newcommand {\Chi} {{\bf \raise 2pt \hbox{$\chi$}} }

\newcommand {\sgn} { {\rm sgn} }

\newcommand {\f}   {\frac}
\newcommand {\p}   {\partial}

\newcommand {\proof} {\noindent {\bf Proof}. }
\newcommand{\beq}{\begin{equation}}
\newcommand{\eeq}{\end{equation}}
\newcommand{\bal}{\begin{align}}
\newcommand{\bc}{\begin{cases}}
\newcommand{\ec}{\end{cases}}
\newcommand{\bea} {\begin{array}{rl}}
\newcommand{\eea} {\end{array}}
\newcommand{\bepa}{\left\{ \begin{array}{l}}
\newcommand{\eepa} {\end{array}\right.}
\newtheorem{theorem}{Theorem}[section]
\newtheorem{lemma}[theorem]{Lemma}
\newtheorem{example}[theorem]{Example}
\newtheorem{definition}[theorem]{Definition}
\newtheorem{remark}[theorem]{Remark}

\newcommand{\qed}{{ \hfill
                     {\unskip\kern 6pt\penalty 500 \raise -2pt\hbox{\vrule\vbox to 6pt{\hrule width 6pt
                     \vfill\hrule}\vrule} \par}   }}
\title{\Large \bf Singular limits for reaction-diffusion equations  with fractional Laplacian and  local or nonlocal nonlinearity}

\author{
Sylvie M\'el\'eard\thanks{CMAP, Ecole Polytechnique, UMR 7641, route de
    Saclay, 91128 Palaiseau Cedex-France; E-mail: \texttt{sylvie.meleard@polytechnique.edu}}, \and Sepideh Mirrahimi\thanks{Institut de Math\'ematiques de Toulouse; UMR 5219, Universit\'e de Toulouse; CNRS, UPS IMT, F-31062 Toulouse Cedex 9, France; E-mail: \texttt{Sepideh.Mirrahimi@math.univ-toulouse.fr}}}

\date{\today}

\begin{document}
\maketitle
\pagestyle{plain}
\pagenumbering{arabic}

\begin{abstract}
\noindent
We perform an asymptotic analysis of models of population dynamics with a fractional Laplacian and local or nonlocal reaction terms. 
The first part of the paper is devoted to the long time/long range rescaling of the fractional Fisher-KPP equation. This rescaling is based on   the exponential speed of propagation of the population. In particular we show that the only role of the fractional Laplacian in determining this speed is at the initial layer where it determines the thickness of the tails of the solutions.\\
Next, we show that such rescaling is also possible for models with non-local reaction terms, as selection-mutation models. However, to obtain a more relevant qualitative behavior for this second case, we introduce, in the second part of the paper, a second rescaling where we assume that the diffusion steps are small. In this way, using a WKB ansatz, we obtain a Hamilton-Jacobi equation in the limit which describes the asymptotic dynamics of the solutions, similarly to the case of selection-mutation models with a classical Laplace term or an integral kernel with thin tails. However, the rescaling introduced here is very different from the latter cases. We extend these results to the multidimensional case. 
\end{abstract}

\section{Introduction}
\label{sec:intro}

We study the asymptotic behavior of the solution of the following equation 
\beq
\label{SM}
\begin{cases}
\p_t n  +  (-\Delta)^{\al/2} n = n\, R(n,I),\\
n(x,0)=n^0(x),\ x\in \mathbb{R}
\end{cases}
\eeq
with 
\beq
\label{def:I}
I(t)=\int_{\mathbb{R}} n(x,t)dx.
\eeq
In all what follows,  $\alpha\in (0,2)$ is given. The term $(-\Delta)^{\al/2}$ denotes the fractional Laplacian:
$$
(-\Delta)^{\al/2} n(x,t) = -\int_0^\infty \left[ n(x+h,t) + n(x-h,t) - 2 n(x,t) \right] \f{dh}{|h|^{1+\al}}.
$$
In the case of the classical diffusion, singular limits using Hamilton-Jacobi equations have been helpful to describe the asymptotic behavior of the reaction-diffusion equations with local or nonlocal nonlinearities (see for instance Freidlin \cite{MFb:85,MF:85}, Evans and Souganidis \cite{LE.PS:89},  Barles et al. \cite{GB.LE.PS:90, GB.SM.BP:09}, Diekmann et al. \cite{OD.PJ.SM.BP:05}). Is it possible to extend these results to the case of models where the Laplace term is replaced by a fractional Laplacian?  \\

\noindent 
Here, we consider  two different forms of reaction terms. The first case corresponds to the fractional Fisher-KPP equation, describing population dispersion with local interactions:
\beq
\label{re-2}
R(n,I)=1-n.
\eeq

\noindent The second case which is motivated by selection-mutation models or spatially structured population models, considers only a dependence on the nonlocal term and therefore induces a nonlocal nonlinearity:
\beq
\label{re-1}
R(n,I)=R(I).
\eeq
A standard example is the logistic one, where $R(I) = (r-I)$, $r$ being the intrinsic growth rate of a population and $I$ is a mean-field competition term. 
Such models are rigorously derived from microscopic (individual-based) dynamics involving L\'evy flights which naturally appear in evolutionary ecology or population dynamics, when the mutation distribution or the dispersal kernel have heavy tails and belong to the domain of attraction of a stable law (see for example Gurney and Nisbet \cite{Gurney} and Baeumer et al. \cite{Baeumer}). This derivation is detailed in Jourdain et al. \cite{Jourdain} for a selection-mutation model where  $x$ denotes a quantitative genetic  parameter. It leads to 
 a  more relevant reaction term with also  a dependence on the trait $x$ of the ecological parameters.  However, due to technical difficulties in this paper we only focus on the nonlocal nonlinearity and study the above simplified version.\\

\bigskip

\noindent 
In this paper, our motivation is twofold. In the one hand, we are interested in the long range/long time asymptotic analysis of  \fer{SM}. The objective is to describe how fast the population propagates, using an asymptotic analysis,  similarly to the case of the classical Fisher KPP equation (see for instance \cite{LE.PS:89,GB.LE.PS:90}). In the other hand, we would like to 
describe the population dynamics, while the mutation (dispersion) steps are small. 
To this end, we look for a rescaling of time and the size of mutation (dispersion)  steps such that we can perform an asymptotic analysis to describe the asymptotic dynamics of the population, similarly to the models where the mutations have thiner tails (see for instance \cite{GB.BP:08,GB.SM.BP:09}). Note that, in the both above cases, usual rescalings used for similar models with the classical laplacian or smoother mutation kernels (see \cite{LE.PS:89,GB.LE.PS:90,GB.BP:08,GB.SM.BP:09}) cannot be used. The possibility of big jumps with a high rate (algebraically small and not exponentially), modifies drastically the behavior of the solutions (see Section \ref{sec:heur}). The leading effect of the big jumps for processes driven by stable L\'evy processes has also been observed from a probabilistic point of view, (see \cite{Imkeller}). 
\\

\noindent 
An important contribution of this paper, is the unusual rescalings that we introduce, which lead to the description of the asymptotic dynamics of the population. 
We will use two rescalings. The first one, being a long range/long time rescaling,  is based on the speed of the propagation.  The idea here, is to look at the population from far, as in homogenization, so that
we forget the full and detailed behavior, but  capture the propagation of the population. In this way, we suggest a new asymptotic formalism to deal with reaction-diffusion equations with the fractional laplacian. This formalism generalizes  a classical approach known as the "approximation of geometric optics" which is well developed in the case of the  reaction-diffusion equations with the classical laplacian (see \cite{MFb:85,MF:85,LE.PS:89,GB.LE.PS:90}).
However, the difference between the speed of propagation in the classical or in the fractional Laplacian case can convince the reader that a similar spatial scaling as in the classical Laplacian case cannot lead to satisfying asymptotical results. Therefore  we introduce a first rescaling  inspired from  the exponential speed of propagation (the result announced in \cite{JR.XC:09} and proved in  \cite{XC.JR:13} by Cabr\'e et al.). 
We show in particular that with this rescaling, the fractional Laplacian disappears in the limit as $\e \to 0$. As a consequence, the only influence of the  fractional Laplacian  on the speed of the propagation is at the initial layer where it determines the thickness of the tails.  This property is very different from the case with the classical Laplacian and is true either in case  \eqref{re-2} or in case \eqref{re-1}. 

\medskip \noindent 
The above rescaling is  not relevant in the case of selection-mutation models or in dispersion models, where we consider small diffusion steps independently of the position of the individuals. Therefore, we suggest a second rescaling for the case of small mutation steps where the mutation (dispersion) kernel is rescaled homogeneously with respect to $x$. In this case the fractional Laplacian does not disappear in the limit as $\e\to 0$ and  the asymptotic dynamics is still influenced by this term.
The asymptotic behavior of the population is described by a Hamilton-Jacobi equation. This approach is closely related to recent works on the asymptotic study of selection-mutation models (see \cite{OD.PJ.SM.BP:05,GB.BP:08,GB.SM.BP:09}) developed in the easier   case where the  mutation steps have finite moments. 
\noindent
For the sake of simple representation, we have presented our results for $x\in \R$. However, we show that these results can be easily extended to the multidimensional cases.
\\



\section{The main results}
\label{sec:heur}

We introduce two scalings yielding two different asymptotics. The first one is a long range/long time rescaling well suited when the equation models a spatial propagation.  The second one is well suited in the selection-mutation modeling, when the diffusion term represents a small mutation approximation. 

\subsection{Long range/long time rescaling and the asymptotic speed of propagation}
\label{sec:heur1}
In this section we  firstly study the asymptotic behavior of   the Fisher-KPP equation \fer{SM} with \fer{re-2}:
\beq
\begin{cases}
\p_t n (x,t) +  (-\Delta)^{\al/2} n = n\, (1-n),\\
n(x,0)=n^0(x),\ x\in \mathbb{R}. \nonumber
\end{cases}
\eeq 
It has been proved in Cabr\'e et al.  \cite{XC.AC.JR:12,XC.JR:13} that the level sets of $n$ propagate with a speed that is exponential in time   (see also  \cite{HE:10,HB.JR.LR:11,JG:11,AC.JR:12} for related works on the speed of propagation {for reaction-diffusion equations, with fractional laplacian or  a diffusion term with thick tails}). In particular, in \cite{XC.JR:13} it is proved that  for any initial data such that
$$
0\leq n_0(x)\leq C\f{1}{1+|x|^{1+\al}},
$$
we have
\beq
\label{roq}
\begin{cases}
n(x,t)\to 0,& \text{uniformly in $\{|x|\geq e^{\sigma t}\}$, if $\sigma >\f{1}{1+\al}$, as $t\to \infty$},\\
n(x,t)\to 1,& \text{uniformly in $\{|x|\leq e^{\sigma t}\}$, if $\sigma <\f{1}{1+\al}$, as $t\to \infty$}.
\end{cases}
\eeq

\noindent Our objective is to understand this behavior using singular limits as for the KPP equation with a Laplace term (cf. \cite{MFb:85,GB.LE.PS:90}). The idea is to rescale the equation and to perform an asymptotic limit so that we forget the full and detailed behavior and capture only this propagation. In the case of the classical Fisher-KPP equation, to study the asymptotic behavior of the solutions one should use the following rescaling \cite{MFb:85,MF:85,LE.PS:89,GB.LE.PS:90}
$$
|x|\mapsto {\f{x}{\e}}, \qquad t\mapsto \f{t}{\e},\quad \text{and}\quad 
n_\e(x,t) = n(  \f{x}{\e} , \f{t}{\e} ).
$$

\noindent

\noindent In the case of the fractional Fisher-KPP equation,  being inspired from \fer{roq} we use the following  long-range/long time rescaling
\beq
\label{res}
|x|\mapsto |x|^{\f{1}{\e}}, \qquad t\mapsto \f{t}{\e}.
\eeq

\noindent
For  the sake of simple representation we assume
\beq
\label{as:sym}
n_0(x)=n_0(|x|),\qquad 
x\in \R.
\eeq
Having in mind that under assumption \fer{as:sym}, for all $(x,t) \in \R \times \R^+$, we have $n(x,t)=n(|x|,t)$, we then can define
\beq
\label{res-n}
n_\e(x,t) = n(  |x|^{\f{1}{\e}} , \f{t}{\e} ).
\eeq
Note that Assumption \fer{as:sym} is not necessary for our results to be held. In the case where $n$ is not symmetric, it is enough to perform the following rescaling
$$
n_\e(x,t) = n(  \sgn (x) \,|x|^{\f{1}{\e}} , \f{t}{\e} ).
$$

\noindent
Replacing \fer{res-n} in \fer{SM} we obtain, 
\beq
\label{KPPe}
\begin{cases}
\e \p_t n_\e (x,t) = \int_0^\infty \left( n_\e \left( \left| |x|^{\f 1\e}+h \right|^\e,t \right) + n_\e \left( \left| |x|^{\f 1\e}-h \right|^\e,t \right)  -2n_\e(x,t) \right) \f{dh}{|h|^{1+\al}}+n_\e(x,t) (1 - n_\e(x,t)),\\
n_\e(x,0)=n_\e^0(x),
\end{cases}
\eeq
where $I_\e(t)=I(\f t \e)$. 

\medskip \noindent 
Although, for the classical Fisher-KPP equation, the long range and long time rescaling coincides with the one with small diffusion steps and long time, this is not the case for  the fractional Fisher-KPP equation. To understand this better, we rewrite \fer{KPPe}, for $x\neq 0$, in the following form
\begin{eqnarray*}
\e \p_t n_\e (x,t) &=& |x|^\f{-\al}{\e}\,\int_0^\infty \left( n_\e \left(|x|\cdot e^{\e k},t \right) + n_\e \left(|x|\cdot \exp\left( \e \log |2- e^{k} | \right),t \right)  -2n_\e(x,t)\right) \f{ e^{k }}{|e^{k } -1|^{1+\al}}dk \\
&&+ \,  n_\e(x,t)\,  (1 - n_\e(x,t)).
\end{eqnarray*}
Notice that here we have used the following change of variable:
$$
h=|x|^\f{1}{\e}\left(e^{k} -1\right),
\quad
\text{so that} 
\quad
|x|^{\f{1}{\e}}+h=|x|^{\f{1}{\e}} \cdot e^{k}.
$$
On this form, one can guess that the fractional Laplatian will disappear in the limit. Note that, by a change of variable, this rescaling can be interpreted as a rescaling of  the integral kernel:
$$
h=|x|(e^k-1),\quad  t\mapsto \f{t}{\e},\quad    M(x,k,dk)=\f{|x|^{-\al}e^{k } dk}{|e^{k }-1|^{1+\al}} \mapsto M_\e(x,k,dk)= |x|^{-\f \al \e}\f{e^{\f k \e} \f{dk}{\e}}{|e^{\f k\e }-1|^{1+\al}}.
$$
We observe that this rescaling is heterogeneous in $x$, and the diffusion steps are rescaled differently at different points. \\

\noindent
 Another way to have an idea of the shape of the solutions is to recall the following bounds
on the   transition probability function $p$ associated with the fractional Laplacian with coefficient $\alpha/2$ (see, e.g., Sato \cite{Sato} p.89 and p.202): 
\beq
\label{p}
\f{ B_m }{t^{\f{1}{\al}}(1+|t^{\f{-1}{\al}}x|^{1+\al})}\leq p(x,t) \leq \f{ B_M }{t^{\f{1}{\al}}(1+|t^{\f{-1}{\al}}x|^{1+\al})}.
\eeq
Note that  the solution $v$ to the following equation
$$
\begin{cases}
\p_t v +(-\Delta)^{\al/2} v= 0,& \text{in $\R \times \R^+$,}\\
v(x,0)=v^0,& \text{in $\R$,}
\end{cases}
$$ 
satisfies
$$
v(x,t)=\int_\R p(y,t) v^0(x-y)dy.
$$
The inequality \fer{p}  is written after the rescaling \fer{res} as
\beq
\label{pe}
\f{ B_m }{(\f{t}{\e})^{\f{1}{\al}}(1+|(\f{t}{\e})^{\f{-1}{\al}}|x|^{\f{1}{\e}}|^{1+\al})}
\leq 
p_\e(|x|,t)=p(|x|^{\f{1}{\e}},\f t \e) 
\leq 
\f{ B_M }{(\f{t}{\e})^{\f{1}{\al}}(1+|(\f{t}{\e})^{\f{-1}{\al}}|x|^{\f{1}{\e}}|^{1+\al})},
\eeq

\medskip\noindent
Being inspired now by \fer{pe} we use the  classical Hopf-Cole transformation

\beq
\label{Hopf}
n_\e=\exp\left( \f{u_\e}{\e} \right),
\eeq
and make the following assumption
\beq
\label{as:inikpp}
\f{C_m}{1+|x|^{\f{1+\al}{\e}}} \leq n_\e(x,0)\leq \f{C_M}{1+|x|^{\f{1+\al}{\e}}}, \quad \text{with $C_m<1<C_M$}.
\eeq

\bigskip
\noindent
Our first result is the following. 

\begin{theorem}\label{th:kpp}
Let $n_\e$ be the solution of \fer{KPPe} with \fer{re-2} and $u_\e=\e\log n_\e$. (i) Under assumption \fer{as:sym} and \fer{as:inikpp}, as $\e\to 0$, $(u_\e)_\e$ converges locally uniformly to $u$ defined as below
\beq
\label{udef}
u(x,t)=\min( 0, -(1+\al) \log |x|+t).
\eeq
(ii) Moreover, as $\e\to 0$,
\beq
\label{limn}
\begin{cases}
n_\e\to 0, & \text{locally uniformly in  $\mathcal A=\{(x,t)\in \R\times (0,\infty) \, |\, t <(1+\al)\log |x|\}$,}\\
 n_\e \to 1, & \text{locally uniformly in  $\mathcal B=\{(x,t)\in \R\times (0,\infty) \, |\, t >(1+\al)\log |x|\}$.}
\end{cases}
\eeq
\end{theorem}

\bigskip \noindent Let us provide some heuristic arguments to understand this result. Rewriting  \fer{KPPe} in terms of $u_\e$ we find
$$
 \p_t u_\e (x,t) =\int_0^\infty \left( \dfrac{n_\e \left( \left| |x|^{\f 1\e}+h \right|^\e,t \right) }{n_\e(x,t) }+ \dfrac{n_\e \left( \left| |x|^{\f 1\e}-h \right|^\e,t \right)}{n_\e(x,t)}  -2 \right) \f{dh}{|h|^{1+\al}}+1-n_\e.
 $$
 To prove the convergence of $(u_\e)_\e$ in Theorems \ref{th:kpp} and \ref{th:RI}, a key point is to find appropriate sub and supersolutions for \fer{KPPe}.
 In Section \ref{sec:pr1} we will prove that the first term of the r.h.s. of the above equation vanishes as $\e\to 0$. 
Therefore, the only remaining term is the one coming from the reaction term, i.e. $u$ the limit of $(u_\e)_\e$, satisfies
\beq
\label{usimple}
\max (\p_t u-1,u)=0.
\eeq
Note that, here the variational form of the equation comes from the fact that $n_\e$ is bounded. 
Suppose now that the initial data in \fer{SM} satisfies
\beq
\label{inifrac}
\f{C_m}{1+|x|^{\f{1+\al}{\e}}}\leq n_\e^0(x) \leq \f{C_M}{1+|x|^{\f{1+\al}{\e}}}.
\eeq
Indeed, note in view of \fer{p} that starting with a compactly supported initial data $n_\e^0$, the tails of  $n_\e$ would have algebraic tails as above, for all $t>0$. \\
Combining the above inequalities, we find  that
$$
 u(x,t)=\min \left(0,-(1+\al) \log |x|+t\right).
$$
The above equality, and the fact that the only steady states of the reaction term \fer{re-2} are $0$ and $1$, suggest that
$$
n_\e \to \begin{cases}
1& \text{in  $\ \mathcal A$,}\\
0& \text{in  $\ \mathcal B$,}
\end{cases}
$$
which is in accordance with \fer{roq}.\\\noindent
 To prove the convergence of $n_\e$ in Theorems \ref{th:kpp} the difficulty is  for the set $\mathcal{B}$. To prove the convergence of $n_\e$ in this set, being inspired by the results on the classical Fisher-KPP equation (see \cite{LE.PS:89}), we introduce an appropriate viscosity (supersolution) test function which leads to the result.

\medskip \noindent
In view of \fer{usimple}, we notice that at the limit $\e=0$, the fractional Laplacian does not have any impact on the dynamics of $u$ and the dynamics are determined only by the reaction term. The only role of the fractional Laplacian in the limit is at the first initial time where the tail of the solution is forced to satisfy  some inequalities similar to \fer{inifrac}. The exponential propagation is hence derived only from the form of the solution at the initial layer. This is an important difference with the KPP equation with the classical Laplacian where, the Laplace term not only forces the solution to have an exponential tail but also it still influences in positive times the dynamics and modifies the speed of propagation. To observe this property consider the following equation
$$ 
\begin{cases}
\p_t m-\da \Delta m=m(1-m),\qquad \da \in \{0,1\},\\
m(x,0)=\exp\left(-\f{x^2}{2}\right).
\end{cases}
$$
It is easy to verify that in long time the invasion front scales as $x \sim \sqrt{2t}$ for $\da=0$, while  for $\da=1$ the invasion front scales as $x\sim 2t$. Therefore, the diffusion term speeds up the propagation. \\
Next we consider an analogous equation but with fractional Laplacian:
$$ 
\begin{cases}
\p_t m+\da (- \Delta)^{\f \al  2} m=m(1-m),\qquad \da \in \{0,1\}\quad \al\in [0,2],\\
m(x,0)=m^0(x),\quad \text{with $m^0$ satisfying \fer{inifrac}}.
\end{cases}
$$
Then following the computations above, in long time and for both cases $\da=0,1$, the invasion front scales as $x \sim e^{\f{t}{1+\al}}$.

\bigskip
\noindent
Let us now state that such  rescaling is also possible for Equation \eqref{SM} with nonlocal interactions \fer{re-1}:
\beq
\begin{cases}
\p_t n  +  (-\Delta)^{\al/2} n = n\, R(I),\\
n(x,0)=n^0(x),\ x\in \mathbb{R}. \nonumber
\end{cases}
\eeq 
For the case \fer{re-1} we additionally assume
\beq
\label{as:maxR}
R(n,I)=R(I),\qquad R(I_0)=0, \quad \text{for some positive constant $I_0>0$}.
\eeq
\beq
\label{as:monR}
-C_1\leq\f{d}{d I} R(I) \leq -C_2,\qquad \text{for all $I\in \R^+$ and positive constants $C_{1}, C_{2}$,}
\eeq
\beq
\label{as:I}
I_m\leq I_\e(0)\leq I_M, \quad \text{where $I_m$ and $I_M$ are positive constants such that  $I_0\in [I_m,I_M]$.}
\eeq 

\begin{theorem}
\label{th:RI}
Let $n_\e$ be the solution of \fer{KPPe} with \fer{re-1} and $u_\e=\e\log n_\e$. (i) Under assumptions \fer{as:sym}, \fer{as:inikpp},  \fer{as:maxR}, \fer{as:monR} and \fer{as:I}, as $\e\to 0$, $(u_\e)_\e$ converges locally uniformly to $u\in \mathcal{C}(\R)$ defined as below
$$
u(x,t)=\min(0,-(1+\al) \log |x|).
$$
(ii) Moreover, $n_\e$ converges, along subsequences as $\e\to 0$, in $L^\infty$ weak-$\ast$ to a function $n\in L^\infty(\R\times \R^+)$, such that $\mathrm{supp} \; n\subset \{(x,t)\in \R\times \R^+ \, | \, u(x,t)=0\}= [-1,1]\times \R^+ $.
\end{theorem}

\subsection{Diffusion with small steps and long time}
\label{sec:heur2}

\noindent  The rescaling  \eqref{res-n} is not satisfying for the case of structured population dynamics. In that case, $x$ represents an hereditary parameter as a phenotypic trait and  the fractional Laplatian term corresponds to a mutation term where an individual with trait $x$ can give births to individuals with traits $x+h$ or $x-h$. The fractional Laplatian models  large mutation jumps (see Jourdain et al. \cite{Jourdain}). We are interested in the long time behavior of the populations with small mutation steps. On the one hand, this rescaling is not adapted to study a solution $n_\e$ which is close to a Dirac mass, while this is  likely the case in selection-mutation models. The  rescaling \eqref{res-n} attributes indeed to $n_\e(x,t)$, for $x\in (-1,1)$ a value close to $n(0,t)$ and therefore flattens the solution. 
On the second hand, since in the context of selection-mutation models, the rescaling is on the size of the mutations and not on the variable $x$ to consider the long range limit, the non homogeneity in the mutation kernel induced by \eqref{res-n}  is not realistic. 
Therefore, in this case, we consider the following rescaling where the mutation kernel remains independent of $x$:
\begin{equation}
\label{res-2}
h=e^k-1,\quad  t\mapsto \f{t}{\e},\quad    M(k,dk)=\f{e^{k } dk}{|e^{k }-1|^{1+\al}} \mapsto M_\e(k,dk)= \f{e^{\f k \e} \f{dk}{\e}}{|e^{\f k\e }-1|^{1+\al}}.
\end{equation}
In this way the size of mutations is rescaled to be smaller homogeneously and independently of $x$.
Here, we have also made a change of variable in time, to be able to observe the effect of small mutations on the dynamics.
An advantage of this choice is that  the size of mutations are rescaled to be smaller homogeneously and independently of $x$.   In this way, $n_\e$ solves the following equation
\beq
\label{SMe}
\begin{cases}
\e \p_t n_\e (x,t) = \int_0^\infty \left( n_\e(x+e^{\e k}-1,t) + n_\e(x-e^{\e k}+1,t) - 2 n_\e(x,t) \right) \f{e^kdk}{|e^k-1|^{1+\al}} +n_\e(x,t) \,R(n_\e,I_\e)(x,t),\\
n_e(x,0)=n_\e^0(x),
\end{cases}
\eeq
with 
\beq
I_\e(t)=\int_{\mathbb{R}} n_\e(x,t)dx,
\eeq
and
$$
n_\e(x,t)=n(x,\f t \e).
$$
 Note that, the rescaling does not change the algebraic distribution of jumps and the problem remains different from what is studied previously. 
 \\

\noindent
We use the following assumptions
\beq
\label{as:ue0}
(u_\e^0)_\e \text{ is a sequence of Lipschitz continuous functions which converge in $C_{\rm loc}(\R)$ to $u^0$, as $\e \to 0$,}
\eeq
\beq
 \text{there exist positive constants $\, A<\al\,$ and $\,B\,$ such that   $\forall \e>0$ and $\forall x, h\in \R$,}\nonumber\; 
\eeq
\begin{eqnarray}
\label{as:u0}
u_\e^0(x)&\leq& -A \log(|x|+1) +B,\\
\label{as:biz}
u_\e^0(x+h) &\leq& u_\e^0(x)+A \log \left(1+|h|\right).
\end{eqnarray}
Note that the above assumptions are satisfied for instance for $u_\e^0(x)=-A\log(|x|+1)+B.$ 

\bigskip \noindent
Rescaling \eqref{res-2} being motivated by mutation-selection models, we will first focus on the the case \fer{re-1}  with nonlocal nonlinearity as developed in \cite{Jourdain}.

\begin{theorem}\label{th:main}
Let $n_\e$ be the solution of \fer{SMe} with \fer{re-1} and $u_\e=\e\log n_\e$. Assume \fer{as:maxR}, \fer{as:monR}, \fer{as:I}, \fer{as:ue0}, \fer{as:u0} and  \fer{as:biz}. (i) Then, as $\e\to 0$,  $(I_\e)_\e$ converges locally uniformly to $I_0$ and $(u_\e)_\e$ converges locally uniformly to a  continuous function $u$ which is  Lipschitz continuous with respect to $x$ and continuous in $t$. Moreover, $u$ is the unique viscosity solution to the following equation
\beq
\label{HJ}
\begin{cases}
 \p_t u-\int_0^\infty \left( e^{D_x  u\cdot k}+e^{-D_x  u\cdot k}-2 \right) \f{e^kdk}{|e^k-1|^{1+\al}} =0,\\
u(x,0)=u^0(x),
\end{cases}
\eeq
and
\beq
\label{max}
\|D_x u \|_{L^\infty(\R \times \R^+)}\leq A,\qquad \max_{x\in \R} u(x,t)=0.
\eeq
(ii) Finally, along subsequences as $\e\to 0$, $n_\e$ converges in $L^\infty\left( w*(0,\infty) ; \mathcal{M}^1(\R) \right)$ to a measure $n$, such that, $supp\, n \subset \{(x,t) \,|\, u(x,t)=0\}$. 
 \end{theorem}
 
 \noindent
 Here, we observe that, contrarily to the case of long range/long time rescaling (Theorems \ref{th:kpp} and \ref{th:RI}), the fractional Laplacian does not disappear at the limit and still has an influence on the dynamics of the asymptotic solution. Later in Example \ref{ex:1},  we give more details on the signification of the results in Theorem \ref{th:main}, for a particular case.\\
 
 \noindent
The above result can be compared to Theorem 1.2 in \cite{GB.SM.BP:09}, where a similar problem has been studied but with a mutation kernel with exponentially small tails. Although, the results seem closely related, the main difference comes from the rescaling that we have considered to obtain such Hamilton-Jacobi equations. In the case of the fractional Laplacian, we should contract much more the mutation steps to obtain a limiting behavior. Note also that, in \cite{GB.SM.BP:09}, in the Hamilton-Jacobi equation obtained at the limit, there is still a dependency in $I(t)$ which is the limit of $(I_\e(t))_\e$ since in that case, the growth rate depends on $x$. 
 \\
 
 \begin{remark}The result in Theorem \ref{th:main}-(ii) can be improved. One can indeed use arguments similar to the one in \cite{GB.BP:08}(Section 3) and the fact that $R(I_\e)$ is small, to obtain that $n_\e$ converges, along subsequences as $\e\to 0$, in 
 $\mathcal{C}\left( (0,\infty) ; \mathcal{M}^1(\R) \right)$ to a measure $n$, and hence for all $t>0$, $supp\, n(\cdot,t) \subset \{u(\cdot,t)=0\}$. However, in this paper, we do not give the proof of this stronger result since we want to focus on the difficulties coming from the nonlocal diffusion. 
 \end{remark}

\noindent Let us now study the case  \fer{re-2}.  In addition to the previous assumptions, we also assume that there exists a positive constant $C_M$ such that
\beq
\label{as:n01}
0\leq  n_\e(x,0)\leq C_M, \qquad \text{for all $x\in \R$ and $\e>0$.}
\eeq

\medskip
\begin{theorem}\label{th:main2}
Let $n_\e$ be the solution of  \fer{SMe} with \fer{re-2}  and $u_\e=\e\log n_\e$. Assume \fer{as:ue0}, \fer{as:u0}, \fer{as:biz} and \fer{as:n01}. (i) Then, as $\e\to 0$, $(u_\e)_\e$ converges locally uniformly to a  function $u$ that is Lipschitz continuous with respect to  $x$ and continuous in $t$. Moreover, $u$ is the viscosity solution to the following Hamilton-Jacobi equation
\beq
\label{HJ2}
\begin{cases}
\max \left( \p_t u-\int_0^\infty \left( e^{D_x  u\cdot k}+e^{-D_x  u\cdot k}-2 \right) \f{e^kdk}{|e^k-1|^{1+\al}}-1 ,u\right)=0,\\
u(x,0)=u^0(x),
\end{cases}
\eeq 
and
\beq
\label{LipKPP}
 \|D_x u \|_{L^\infty(\R \times \R^+)}\leq A.
\eeq
 (ii) Moreover, as $\e\to 0$,
\beq
\label{limn2}
\begin{cases}
n_\e\to 0, & \text{locally uniformly in  $\{(x,t)\in \R\times (0,\infty) \, |\, u(t,x)<0 \}$,}\\
 n_\e \to 1, & \text{locally uniformly in  $\mathrm{Int}\;\{(x,t)\in \R\times (0,\infty) \, |\,u(t,x)=0 \}$.}
\end{cases}
\eeq

\end{theorem}
See Example \ref{ex:2}, for an interpretation of the results in the above Theorem, for a particular case.\\

\noindent
To prove Theorems \ref{th:main} and \ref{th:main2} we first prove some regularity bounds using some sub- and supersolution arguments. Next, to prove the convergence to the corresponding Hamilton-Jacobi equation, we use the so called  half-relaxed limits method for viscosity solutions, see \cite{GB.BP:88}. Note that, Theorems \ref{th:main} and \ref{th:main2} are proved under the thick tail assumption \fer{as:biz} on the initial data, which assumes that $n_\e^0$ has tails of order $|x|^{-A/\varepsilon}$. Indeed, we need such property to be able to pass to the Hamilton-Jacobi limit. We still don't know how the solution would behave in the case where the initial data has a thiner tail. Note that for the Hamiltonian in \fer{HJ} (or in \fer{HJ2}) to be finite, one should at least have $|D_x u|< \al$ in $\R\times (0,\infty)$.\\

\noindent
Note that, our asymptotic study, or more generally the "approximation of geometric optics" approach is closely related to the large deviation theory (see for instance \cite{MFb:85,MF:85}). In \cite{CB.EC:09,CB.EC:13}, some large deviation type estimates have indeed been proven for some nonlocal equations with Levy type kernels which have fast decays. In those papers, the kernel must scale at most as $e^{-|x|}$ and therefore, the case of the fractional Laplacian is not treated. It is however worth mentioning that, with our second rescaling in Theorems   \ref{th:main} and \ref{th:main2}, although at the $\e$ level it is not the case, in the limit $\e=0$, the problem approaches the case of Levy type kernels which scale as $e^{-\al |x|}$ (known as kernels with tempered stable law) and we obtain a Hamilton-Jacobi equation with a Hamiltonian similar to the one obtained in \cite{CB.EC:13}.\\

\bigskip \noindent 
Lest us  provide heuristic arguments on the proof of  Theorem \ref{th:main}.  
Replacing \eqref{res-2} in \fer{SM}  one obtains \fer{SMe}. Then, using the Hopf-Cole transformation \fer{Hopf}, in the case of \fer{re-1} one obtains
$$
 \p_t u_\e (x,t) =  \int_0^\infty \left( e^{\f{u_\e(x+e^{\e k}-1,t)-u_\e(x,t)}{\e}} + e^{\f{u_\e(x-e^{\e k}+1,t)-u_\e(x,t)}{\e}}-2\right) \f{ e^{ k} }{ | e^{ k} -1|^{1+\al}}dk+ R(I_{\e}(t)).
$$
which leads formally to
$$
 \p_t u=\int_0^\infty \left( e^{D_x  u\cdot k}+e^{-D_x  u\cdot k}-2 \right) \f{e^kdk}{|e^k-1|^{1+\al}}+R(I),
$$
where $I$ is the limit of $(I_\e)_\e$ as $\e\to 0$. Our conjecture is that an equivalent result as in Theorem \ref{th:main} holds for the general case $R(x,I)=R(I)$ but  due to technical  difficulties this is beyond the scope of the present paper.\\

\noindent In the case of \fer{re-2} following similar computations as above we find formally
$$
\max \left( \p_t u-\int_0^\infty \left( e^{D_x  u\cdot k}+e^{-D_x  u\cdot k}-2 \right) \f{e^kdk}{|e^k-1|^{1+\al}}-1 ,u\right)=0.
$$

\bigskip
\noindent Next, we give two examples where we discuss the interpretation of the results in Theorems \ref{th:main} and \ref{th:main2}.

\begin{example}
\label{ex:1}
Let $n_\e$ be the solution of \fer{SMe} with \fer{re-1} and $n_\e^0(x)= \left( |x|+1 \right)^{-A/\e }$. Then, it follows from Theorem \ref{th:main} that $(u_\e)_\e$ converges locally uniformly to the unique viscosity solution of the following equation
$$
\begin{cases}
\p_t u -H(D_x u)=0,\\
u(x,0)=-A \log \left( 1+ |x| \right),
\end{cases}
$$
with 
$$
H(D_x u)= \int_0^\infty \left( e^{D_x u\cdot k}+ e^{-D_x u\cdot k} -2 \right) \f{e^k dk}{|e^k-1|^{1+\al}}.
$$
It is easy to verify, using a Taylor expansion, that
$$
\underline{C}\, p^2=p^2\, \int_0^\infty \f{k^2 \left( e^{-Ak}+1\right)e^k}{2|e^k-1|^{1+\al}}dk
\leq  H(p) \leq p^2 \,\int_0^\infty \f{k^2 \left( e^{Ak}+1\right)e^k}{2|e^k-1|^{1+\al}}dk=\overline{C}\, p^2.
$$
The above estimates, allow us to approximate the value of $u$:
$$
 \sup_{y\in \R} \, \{-A \log \left(1+|y| \right) -\f{|x-y|^2}{4\underline{C}t} \} \leq u(x,t) \leq \sup_{y\in \R} \, \{-A \log \left(1+|y| \right) -\f{|x-y|^2}{4\overline{C}t} \}.
$$
In particular, it follows that 
$$
\mathrm{supp}\; n =\{0\}\times \R^+,\quad \text{and thus for all $(x,t)\in \R\times \R^+$,}\quad n(x,t)=I_0\,\delta(x=0).
$$
Moreover, at the other points ($x\neq 0$), $u$ becomes more and more flat as time goes by, and finally as $t\to \infty$, $u(x,t)\to 0$ for all $x\in \R$.

\noindent We note that here $n$ does not evolve and it remains a Dirac mass at $0$, since the growth rate $R$ is too simple. With the present form of $R$ there is no reason for the population to move from one point to another. For the Dirac mass to evolve in time, the growth rate $R$ must depend on $x$, as  was the case for instance in \cite{GB.BP:08,GB.SM.BP:09}.
\end{example}

\begin{example}
\label{ex:2}
Let $n_\e$ be the solution of \fer{SMe} with \fer{re-2} and $n_\e^0(x)= \left( |x|+1 \right)^{-A/\e}$. Then, it follows from Theorem \ref{th:main2} that $(u_\e)_\e$ converges locally uniformly to the unique viscosity solution of the following equation
$$
\begin{cases}
\max \left( \p_t u -H(D_x u)+1,u \right)=0,\\
u(x,0)=-A \log \left( 1+ |x| \right),
\end{cases}
$$
with $H(D_x u)$ defined in Example \ref{ex:1}. Using the estimates presented in Example \ref{ex:1}, we obtain that
$$
 \min \left( \sup_{y\in \R} \, \{-A \log \left(1+|y| \right) -\f{|x-y|^2}{4\underline{C}t}+t \},0 \right) \leq u(x,t) \leq \,  \min \left( \sup_{y\in \R} \{-A \log \left(1+|y| \right) -\f{|x-y|^2}{4\overline{C}t} +t\} ,0 \right).
$$
After some computations, we find
$$
\begin{array}{c}
\left\{(x,t)\in \R\times \R^+\, |\,
|x|\leq \max_{r\in [0,1]} \left[2\sqrt{\underline{C}}rt+e^{\f{t(1-r^2)}{A}}-1\right] \right\}
\subset
\{u=0\}\\
\subset \left\{(x,t)\in \R\times \R^+\, |\,
|x|\leq \max_{r\in [0,1]} \left[2\sqrt{\overline{C}}rt+e^{\f{t(1-r^2)}{A}}-1\right] \right\}.
\end{array}
$$
In view of \fer{limn2}, the above line indicates that the population propagates in space and that the front position still moves exponentially in time.

\end{example}

\bigskip
\noindent
The remaining part of the article is organized as follows.  In Section \ref{sec:prel} we give some preliminary results on the boundedness of $n_\e$ and $I_\e$. Section \ref{sec:pr1} is devoted to the proofs of  Theorems \ref{th:kpp} and \ref{th:RI}. In Section \ref{sec:reg1} we prove some regularity results for \fer{SMe} with the reaction term given by \fer{re-1}. In Section \ref{sec:reg2} we prove some regularity results for \fer{SMe} with the reaction term given by \fer{re-2}. We next prove Theorems \ref{th:main} and \ref{th:main2} respectively in Sections \ref{sec:th1} and \ref{sec:th2}. Finally, we show how our results can be extended to the multidimensional case in Section \ref{sec:multid}.

\medskip \noindent 
Throughout the paper, we denote by $C$ positive constants that are independent of $\e$ but can change from line to line. The notion of solutions that we consider throughout the paper, is classical unless stated otherwise.
\section{Notations and preliminary results }
\label{sec:prel}

It is classical that \fer{SM} with \fer{re-2} has a unique classical solution which is smooth. The existence of a unique weak solution for  a more general equation than \fer{SM} with \fer{re-1}, is  proved in \cite{Jourdain}. Moreover, from the regularizing effect of the fractional Laplacian we deduce that the solution is smooth and hence classical. We prove additionally some uniform bounds on $n_\e$ and $I_\e$ respectively in cases \fer{re-2} and \fer{re-1}, which are derived from the monotonicity in the reaction term.

\begin{lemma}\label{th:n}
Let $n_\e$ be the unique solution of  \fer{KPPe} or \fer{SMe} with \fer{re-2}.
Under assumption \fer{as:n01}, we have 
\beq
\label{boundn}
0 \leq n_\e(x,t) \leq \f{C_M e^{\f t \e}}{1-C_M+C_Me^{\f t \e}},\qquad \text{for all $(x,t)\in \R\times [ 0,\infty)$}.
\eeq
\end{lemma}

\proof
One can easily verify that  the nul function  is a subsolution and the r.h.s. is a supersolution of  \fer{KPPe} and \fer{SMe} with \fer{re-2}. Hence, \fer{boundn} follows from Assumption \fer{as:inikpp} or Assumption \fer{as:n01} and the comparison principle. 
\qed

\begin{lemma}\label{th:I}
Let $n_\e$ be the unique solution of  \fer{KPPe} or \fer{SMe} with \fer{re-1}.
Under assumptions  \fer{as:maxR}, \fer{as:monR} and \fer{as:I}, we have 
\beq
\label{boundI}
I_m \leq I_\e(t) \leq I_M,\qquad \text{for all $t\geq 0$}.
\eeq
Moreover as $\e \to 0$,  $(I_\e)_\e$  converges locally uniformly in $(0,\infty)$
 to $I_0$. Moreover, there exists constants $C_3$ and $C_4$ such that
 \beq
 \label{limR}
C_3\,\e \leq  \int_0^t R(I_\e(s))ds \leq C_4\,\e, \qquad \text{for all $t\in \R^+$.}
 \eeq
\end{lemma}

\proof
In both cases of equations \fer{KPPe} and  \fer{SMe} one can obtain
\beq
\label{dtI}
\e\f{d}{dt}I_\e(t)=I_\e(t)R(I_\e(t)).
\eeq 
In the case of \fer{SMe}, this can be derived by integrating \fer{SMe} with respect to $x$. In the case of \fer{KPPe}, we integrate first \fer{SM} with respect to $x$ and then make the change of variable $I_\e(t)=I(\f t \e)$.\\

We notice that, using \fer{as:maxR}--\fer{as:monR},
$$
R(I)<0, \quad\text{for all $I>I_0$} \qquad \text{and } \quad R(I)>0,\quad \text{for all $I<I_0$} .
$$
>From the above inequalities and \fer{as:I} we deduce that
\beq
\label{bIe}
I_m\leq \min \left( I_\e(0),I_0\right) \leq I_\e(t) \leq \max \left( I_\e(0),I_0\right)\leq I_M,\qquad \text{for all $t\geq 0$}.
\eeq
Moreover,  $I_\e(t)$ is monotone in $\R^+$, since $R(I_\e(t))$ does not change sign in this interval. We now suppose that $I_\e(0)<I_0$. The case with $I_\e(0)>I_0$ can be studied following similar arguments.
We compute using \fer{as:monR}
$$
\f{d}{dt} R(I_\e(t))=
\f{d}{dI}R(I_\e(t))\, \f{d}{dt}I_\e(t)\leq -\f{C_2}{\e}I_\e(t)R(I_\e(t)).
$$
Using \fer{bIe}, it follows that 
$$
R(I_\e(t))\leq R(I_\e(0))e^{-\f{C_2 I_m t}{\e}}.
$$
Hence as $\e\to 0$,  $\left( R(I_\e(t)) \right)_\e$ converges  locally uniformly in $(0,\infty)$ to $0$. We then conclude, using again \fer{as:monR} that $\left( I_\e(t) \right)_\e$ converges locally uniformly in $(0,\infty)$ to $I_0$. Moreover, integrating \fer{dtI} we obtain
$$
I_\e(t)=I_\e(0)e^{\f{1}{\e}\int_0^t R(I_\e(s))ds.}
$$
Since $I_\e$ is bounded above and below by positive constants, we obtain \fer{limR}.


\qed

\bigskip


\noindent
To prove our results we will need some comparison principles for  equations of the following type  
\beq
\label{visc}
\p_t n+ r (-\Delta )^{\al/2} n(x,t) +F \left(t,x,n,D_x n\right)=0,\qquad \text{in $\R\times \R^+$,}
\eeq
with $r\geq 0$. We introduce here the statement of the comparison principle that we will use throughout the paper.

\begin{definition}[Comparison principle]
Equation \fer{visc} admits a comparison principle if the following statement holds:\\
Let $n_1$ and $n_2$ be respectively  viscosity subsolution and supersolution of \fer{visc} (see the definition in \cite{GB.CI:08}) and
$$
n_1(x,0) \leq n_2(x,0),\qquad \text{for all $x\in \R$}.
$$
Then 
$$
n_1(x,t) \leq n_2(x,t),\qquad \text{for all $(x,t)\in \R\times \R^+$}.
$$
\end{definition}

\section{The proofs of Theorems \ref{th:kpp} and \ref{th:RI}}
\label{sec:pr1}

\subsection{The proof of Theorem \ref{th:kpp}  }
\label{sec:prkpp}

(i) To prove the first part of Theorem \ref{th:kpp}, we claim that for all $\da >0$, there exists $\e_0(\da)$ small enough such that
\beq
\label{bn}
\f{C_m e^{-\e t-\f{\da}{\e}}}{1+e^{-\f {(t+\da)}{ \e}}|x|^{\f {1+\al}{\e}}} \leq n_\e(x,t)\leq \f{C_M e^{\e t}}{1+e^{-\f {(t+\da)}{ \e}}|x|^{\f {1+\al}{\e}}},\quad \text{for all $\e\leq \e_0$ and in $\R\times \R^+$}.
\eeq
We postpone the proof of the above inequalities to the end of this paragraph.\\

\noindent Combining \fer{Hopf} and \fer{bn} we obtain
\beq
\label{bu}
-\e^2 t - \e \log C_m - \e \log \left(1+e^{-\f{t+\da}{\e}}|x|^{\f{1+\al}{\e}}\right) -\da \leq u_\e (x,t) \leq \e^2 t +\e \log C_M - \e \log \left(1+e^{-\f{t+\da}{\e}}|x|^{\f{1+\al}{\e}}\right).
\eeq
Define
$$
\underline u(x,t)=\liminf_{\e\to 0}u_\e(x,t),\quad \overline u(x,t)=\limsup_{\e\to 0}u_\e(x,t),\quad \text{for all $(x,t)\in \R\times \R^+$}.
$$
Letting $\e\to 0$, we obtain 
$$
\min \left( 0,t+\da-(1+\al) \log |x| \right) -\da \leq \underline u (x,t) \leq \overline u (x,t) \leq \min \left( 0,t+\da-(1+\al) \log |x| \right).
$$
We then let $\da \to 0$ and obtain
$$
 \underline u (x,t)= \overline u (x,t) =\min \left( 0,t-(1+\al) \log |x| \right).
$$
In other words $u_\e$ converges to $u$ given by \fer{udef}. We note that this convergence is locally uniform in $\R \times (0,\infty)$, since for any compact set $K\in \R \times (0,\infty)$  one can pass to the limit in the r.h.s. and the l.h.s. of \fer{bu} uniformly in $\e$.
 \\

\noindent It now remains to prove \fer{bn}. We only prove the r.h.s. of \fer{bn}. The l.h.s. can be proved following similar arguments.\\

\noindent To this end, we need the following lemma, which is proved in Appendix \ref{ap:lem}.
\begin{lemma}
\label{lem:ex}
Let $g(x)=\f{1}{1+|x|^{1+\al}}$. Then, there exists a positive constant $C$, independent of $x$, such that
$$
|(-\Delta)^{\f\al 2} g(x)|\leq C g(x).
$$
\end{lemma}

\bigskip

\bigskip \noindent 
We define
$$
f_\e(x,t):=\f{C_M}{1+e^{-\f{t(1+\e^2)+\da}{\e}}|x|^{\f{1+\al}{\e}}}.
$$
We notice that, for $\e$ small enough, $f_\e$ verifies
$$
\begin{cases}
\f{\p}{\p t} f_\e\geq \f{f_\e}{\e}(1+\e^2-f_\e),\\
f_\e(x,0)=\f{C_M}{1+e^{-\f{\da}{\e}}|x|^{\f{1+\al}{\e}}}.
\end{cases}
$$
Moreover, defining 
$$
\Delta_\e^{\f{\al}{2}} f_\e(x,t) := \int_0^\infty \left( f_\e \left( \left| |x|^{\f 1\e}+h \right|^\e,t \right) + f_\e \left( \left| |x|^{\f 1\e}-h \right|^\e,t \right)  -2f_\e(x,t) \right) \f{dh}{|h|^{1+\al}},
$$
we deduce  from Lemma \ref{lem:ex} and a change of variable that 
$|\Delta_\e^{\f{\al}{2}} f_\e(x,t)| \leq Ce^{-\f{\al\left( (1+\e^2)t+\da \right)}{(1+\al)\e}}f_\e(x,t).$
It follows that for $\e\leq \e_0(\da)$ with $\e_0$ small enough,
$$
|\Delta_\e^{\f{\al}{2}} f_\e(x,t)|\leq \e^2 f_\e(x,t).
$$
We deduce that for all $\e\leq \e_0(\da)$, $f_\e$ is a supersolution of \fer{KPPe} with \fer{re-2}. Moreover $f_\e(x,0)\geq n_\e(x,0)$ thanks to \fer{as:inikpp}. We conclude from the comparison principle for  \fer{KPPe} with \fer{re-2} (see \cite{GB.CI:08} Theorem 3, or \cite{XC.JR:13}) that
$$
n_\e(x,t) \leq  \f{C_M}{1+e^{-\f {(1+\e^2)t+\da}{ \e}}|x|^{\f {1+\al}{\e}}}\leq \f{C_M e^{\e t}}{1+e^{-\f {t+\da}{ \e}}|x|^{\f {1+\al}{\e}}}, \qquad \text{for all $\e\leq \e_0(\da)$}.
$$

(ii) We now prove the second part of Theorem \ref{th:kpp}. We first notice using \fer{udef} that, for any compact set $K\subset \mathcal{A}$, there exists a positive constant $a$ such that for all $(x_0,t_0) \in \mathcal A$ we have $u(x_0,t_0)< -a$. It is thus immediate from \fer{Hopf} that $n_\e$ converges uniformly in $K$ to $0$. Next, we study the case $(x_0,t_0)\in K$, $K$ a compact set such that $K\subset  \mathcal{B}$. To this end, we define
$$
\vp(x,t)= \min \left( 0, -(1+\al)\log |x| +t_0-\da \right) - (t-t_0)^2,
$$
where $\da$ is a positive constant chosen small enough such that for all $(y,s)\in K$, $s\geq 2\da$ and such that $(1+\al)\log |x_0| <t_0-\da$. 
It is easy to verify that $u-\vp$ attains a local in $t$ and global in $x$ minimum at $(x_0,t_0)$. Moreover, this minimum is strict in $t$ but not in $x$. We also define 
$$
\vp_\e(x,t) = -\e \log \left( 1+ e^{-\f{t_0-\da}{\e}} |x|^{\f{1+\al}{\e}} \right)-(t-t_0)^2.
$$
One can also verify that $(\vp_\e)_\e$ converges locally uniformly to $\vp$. Since $(u_\e)_\e$ converges also locally uniformly to $u$, we deduce that there exist points $(x_\e,t_\e) \in K$ such that  $u_\e-\vp_\e$ has a  local in $t$ and global in $x$ minimum at $(x_\e,t_\e)$ and such that $t_\e\to t_0$  and $(u_\e-\vp_\e)(x_\e,t_\e) \to 0$ as $\e\to 0$. 

\noindent We then, using \fer{Hopf}, rewrite \fer{KPPe} as follows
$$
 \p_t u_\e (x,t) =\displaystyle \int_0^\infty \left( e^{\f{u_\e \left( \left| |x|^{\f 1\e}+h \right|^\e,t \right)- u_\e(x,t)}{\e} }+ e^{\f{u_\e \left( \left| |x|^{\f 1\e} - h \right|^\e,t \right)- u_\e(x,t)}{\e} }  -2 \right) \f{dh}{|h|^{1+\al}}+1-n_\e(x,t).
 $$
 Since $u_\e -\vp_\e$ has a  local in $t$ and global in $x$ minimum at $(x_\e, t_\e)$, we have
 $$
  \p_t u_\e (x_\e,t_\e)= \p_t \vp_\e (x_\e,t_\e)=-2(t_\e-t_0),
 $$
$$
\begin{array}{c}
\displaystyle\int_0^\infty \left( e^{\f{u_\e \left( \left| |x_\e|^{\f 1\e}+h \right|^\e,t_\e \right)- u_\e(x_\e,t_\e)}{\e} }+ e^{\f{u_\e \left( \left| |x_\e|^{\f 1\e} - h \right|^\e,t_\e \right)- u_\e(x_\e,t_\e)}{\e} }  -2 \right) \f{dh}{|h|^{1+\al}}\\
\geq
\displaystyle\int_0^\infty \left( e^{\f{\vp_\e \left( \left| |x_\e|^{\f 1\e}+h \right|^\e,t_\e \right)- \vp_\e(x_\e,t_\e)}{\e} }+ e^{\f{\vp_\e \left( \left| |x_\e|^{\f 1\e} - h \right|^\e,t_\e \right)- \vp_\e(x_\e,t_\e)}{\e} }  -2 \right) \f{dh}{|h|^{1+\al}}
=\f{\Delta_\e^{\f{\al}{2}} g_\e(x_\e,t_\e)}{g_\e(x_\e,t_\e)},
 \end{array}
 $$
where $\Delta_\e^{\f{\al}{2}}$ is defined in part (i) of the ongoing proof and 
 $$
 g_\e(x,t) : = \exp \left( \f{\vp_\e(x,t)}{\e} \right)=\f{e^{\f{-(t-t_0)^2}{\e}}} { 1+ e^{-\f{t_0-\da}{\e}}|{x}|^{\f{1+\al}{\e}} }.
 $$
 Using Lemma \ref{lem:ex} and a change of variable similarly to the proof of part (i) we obtain
 $$
\f{ | \Delta_\e^{\f{\al}{2}} g_\e(x_\e,t_\e) |}{g_\e(x_\e,t_\e)} \leq  C e^{-\f{\al(t_0-\da)} {(\al+1)\e}}\leq Ce^{-\f{\al \da} {(\al+1)\e}},
 $$
 which vanishes uniformly for all $(x_\e,t_\e)\in K$ as $\e\to 0$. Combining the above arguments we deduce that
 $\
 n_\e(x_\e,t_\e) \geq 1 +o(1).
 $
 Next, we notice that since $u_\e-\vp_\e$ has a local minimum in $(x_\e,t_\e)$, it follows that
 $$
 u_\e (x_\e,t_\e) -\vp(x_\e,t_\e)\leq u_\e (x_0,t_0) - \vp_\e (x_0,t_0).
 $$
 Moreover, by definition, we have
 $$
 \vp_\e (x_\e,t_\e) \leq \vp_\e(x_0,t_0).
 $$
 Combining the above inequalities we find
 $$
  u_\e (x_\e,t_\e) \leq u_\e(x_0,t_0), \quad \text{and thus }  \quad  n_\e (x_\e,t_\e) \leq n_\e(x_0,t_0).
 $$
We deduce that 
  $$
 n_\e(x_0,t_0) \geq 1 +o(1)  \quad \text{and hence }  \quad \liminf_{\e \to 0} n_\e(x_0,t_0) \geq 1,\quad \text{uniformly in $K$}.
 $$
 Finally, we conclude from the above inequality and Lemma \ref{th:n} that $n_\e(x_0,t_0)\to 1$ uniformly in $K$, as $\e\to 0$. 
 
\subsection{The proof of Theorem \ref{th:RI}  }

\hskip 0.3cm 
(i) The proof of Theorem \ref{th:RI}-(i), is close to the one of Theorem \ref{th:kpp}, (i). In this case, we prove that 
 for all $\da >0$, there exists $\e_0(\da)$ small enough such that
\beq
\label{bnI}
\f{C_m e^{-\e t-\f{\da}{\e}+\f{1}{\e}\int_0^t R(I_\e(s))ds}}{1+e^{-\f {\da}{ \e}}|x|^{\f {1+\al}{\e}}} \leq n_\e(x,t)\leq
 \f{C_M e^{\e t+\f{1}{\e}\int_0^t R(I_\e(s)) ds}}{1+e^{-\f {\da}{ \e}}|x|^{\f {1+\al}{\e}}},\quad \text{for all $\e\leq \e_0$ and in $\R\times \R^+$}.
\eeq
We notice that, admitting the above inequality is true, following similar arguments as in Subsection \ref{sec:prkpp} and using \fer{limR}, we deduce that 
as $\e\to 0$, $(u_\e)_\e$ converges locally uniformly in $\R \times (0,\infty)$ to $u$ defined as below
$$
u(x,t)=u(x)=\min(0,-(1+\al) \log |x|).
$$
\noindent It now remains to prove \fer{bnI}. As before, we only prove the r.h.s. of \fer{bn}. The l.h.s. can be proved following similar arguments. To this end, we define
$$
f_\e(x,t):=\f{C_Me^{\f 1 \e \int_0^t R(I_\e(s))ds+\e t}}{1+e^{-\f{\da}{\e}}|x|^{\f{1+\al}{\e}}}.
$$
We notice that $f$ verifies
$$
\begin{cases}
\f{\p}{\p t} f_\e = \f{f_\e}{\e}\left( R(I_\e)+\e^2 \right),\\
f_\e(x,0)=\f{C_M}{1+e^{-\f{\da}{\e}}|x|^{\f{1+\al}{\e}}}.
\end{cases}
$$
Moreover, we deduce  again from Lemma \ref{lem:ex} that for $\e\leq \e_0(\da)$ with $\e_0$ small enough,
$$
|\Delta_\e^{\f{\al}{2}} f_\e(x,t)| \leq Ce^{-\f{\al\da}{(1+\al)\e}}f_\e(x,t) \leq \e^2 f_\e(x,t).
$$
It follows that for all $\e\leq \e_0(\da)$, $f_\e$ is a supersolution of \fer{KPPe}. Moreover $f_\e(x,0)\geq n_\e(x,0)$ thanks to \fer{as:inikpp}.  We conclude from the comparison principle for  \fer{KPPe} with \fer{re-2} and $I_\e(\cdot)$ fixed (see \cite{GB.CI:08} Theorem 3) that$$
n_\e(x,t)\leq  \f{C_Me^{\f 1 \e \int_0^t R(I_\e(s))ds+\e t}}{1+e^{-\f{\da}{\e}}|x|^{\f{1+\al}{\e}}}, \qquad \text{for all $\e\leq \e_0(\da)$}.
$$

(ii) We first deduce from \fer{bnI} and \fer{limR} that $n_\e$ is uniformly bounded in $L^\infty(\R\times \R^+)$ for all $\e\leq \e_0$. It follows that $n_\e$ converges, along subsequences $\e\to 0$, in $L^\infty$ weak-$\ast$ to a function $n\in L^\infty(\R\times \R^+)$. Moreover, from \fer{Hopf} and the fact that $(u_\e)_\e$ converges locally uniformly to $u$, we deduce that $\mathrm{supp} \; n\subset \{(x,t)\in \R\times \R^+ \, | \, u(x,t)=0)\}=[-1,1]\times \R^+$.

\section{Regularity results for \fer{SMe} and the reaction term given by \fer{re-1}}
\label{sec:reg1}
We first notice that, combining \fer{Hopf} with \fer{SMe}, we obtain
\beq
\label{eq:ue}
\p_t u_\e(x,t) = \int_{0}^\infty \left[ e^{\f{u_\e(x+e^{\e k}-1,t)-u_\e(x,t)}{\e}}+e^{\f{u_\e(x-e^{\e k}+1,t)-u_\e(x,t)}{\e}}-2\right] \f{e^k}{|e^k-1|^{1+\al}}dk+R(I_\e(t)).
\eeq
We then prove the following 
\begin{theorem}
\label{th:reg}
Assume  \fer{as:maxR}, \fer{as:monR}, \fer{as:I},   \fer{as:u0} and \fer{as:biz}. Then, for all $T>0$ and $R>0$, there exist constants $A_1(R,T)$, $A_2(T)$ and $C$ such that
\beq
\label{boundue}
 -\f A 2\log(|x|^2+1)-D-Ct \leq u_\e(x,t)\leq -\f A 2\log(|x|^2+1)+B+Ct,\quad \text{in $ B_R(0)\times [0,T]$},
 \eeq
 and 
 \beq 
 \label{maxue}
 \e \log \left( \f{I_m}{4A_2(T)} \right) \leq \max_{x\in \R} u_\e(x,t),\quad \text{for all $t\in  [0,T]$}.
 \eeq
Moreover, we have
\beq
\label{log-es}
u_\e(x+h,t) \leq u_\e(x,t)+A\log(1+|h|),\qquad \text{for all $x,\, h \in \R$ and $t\geq 0$}.
\eeq
In particular $(u_\e)_\e$ is uniformly Lipschitz with respect to $x$:
\beq
\label{Lipu}
\| D_x u_\e \|_{L^\infty(\R\times \R^+)}\leq A.
\eeq


\end{theorem}

\proof
(i) \textbf{Uniform bound from above.} We prove that, for $C$ large enough,
$$
u_\e(x,t) \leq \overline s(x,t):=-\f A 2\log(|x|^2+1)+B+Ct,\qquad \text{in $ B_R(0)\times [0,T]$}.
$$
We prove indeed that $\overline s$ is a supersolution of \fer{eq:ue}.  One can also verify that,  \fer{eq:ue} with \fer{re-1} and $I_\e$ fixed, admits a comparison principle, since \fer{SMe} admits a comparison principle (see \cite{GB.CI:08} Theorem 3). Then, the claim follows from \fer{as:u0} and since $u_\e$ is a solution and in particular a subsolution of \fer{eq:ue}.\\

\noindent To prove that $\overline s$ is a supersolution of \fer{eq:ue}, since $R$ is bounded thanks to \fer{as:monR} and \fer{boundI}, it is enough to prove that, for $C$ sufficiently large but independent of $\e$,
 $$
 S:=\int_{k\geq 0} \left[ e^{\f{\overline s(x+e^{\e k}-1,t)-\overline s(x,t)}{\e}}+e^{\f{\overline s(x-e^{\e k}+1,t)-\overline s(x,t)}{\e}}-2\right] \f{e^k}{|e^k-1|^{1+\al}}dk\leq C.
$$
We compute
$$\begin{array}{rl}
S&=\displaystyle\int_{k\geq 0} \left[ 
\dfrac { \left( |x|^2+1 \right)^{\f{A}{2\e} }}  {\left( |x+e^{\e k}-1|^2+1 \right)^{\f{A}{2\e} } }+ \dfrac { \left( |x|^2+1 \right)^{\f{A}{2\e} }} {  \left( |x-e^{\e k} + 1|^2+1 \right)^{\f{A}{2\e} }} -2\right] \f{e^k}{|e^k-1|^{1+\al}}dk = f+g,
\end{array}
$$
with
\beq
\label{def:f}
f=\int_{1\geq k\geq 0}
\left[  \dfrac { \left( |x|^2+1 \right)^{\f{A}{2\e} }}  {\left( |x+e^{\e k}-1|^2+1 \right)^{\f{A}{2\e} } }+ \dfrac { \left( |x|^2+1 \right)^{\f{A}{2\e} }} {  \left( |x-e^{\e k} + 1|^2+1 \right)^{\f{A}{2\e} }} -2\right] \f{e^k}{|e^k-1|^{1+\al}}dk,
\eeq
and
\beq
\label{def:g}
g= \displaystyle\int_{k\geq 1} \left[ \left( \dfrac { |x|^2+1 }  { |x+e^{\e k}-1|^2+1 }  \right) ^{\f{A}{2\e} } 
+  \left( \dfrac { |x|^2+1 }  { |x-e^{\e k}+1|^2+1 }  \right) ^{\f{A}{2\e} }  -2\right] \f{e^k}{|e^k-1|^{1+\al}}dk.
\eeq

\noindent Let
$$
s_{\e,x}^1(k)=\left(\dfrac {|x|^2+ 1 } {|x+e^{\e k}-1|^2+1} \right) ^{\f{A}{2\e} },\quad 
s_{\e,x}^2(k)=\left( \dfrac {|x|^2+ 1 } {|x-e^{\e k}+1|^2+1} \right) ^{\f{A}{2\e} }.
$$
We claim that
\beq
\label{sepx}
s_{\e,x}^i(k)\leq e^{Ak},\qquad \text{for $i=1,2$ and all $\e>0$, $k\geq 0$ and $x\in \R$.}
\eeq
We show this only for $i=1$. The case $i=2$ can be proved following similar arguments.\\
For the sake of simple representation we introduce a new variable 
$$
y=x+l, \qquad l=e^{\e k}-1.
$$
Then $s_{\e,x}^i(k)$ is rewritten in terms of $y$ and $l$ as 
$\ 
s_{\e,x}^i(k) = \left(\dfrac {|y-l |^2+ 1 } {|y|^2+1} \right) ^{\f{A}{2\e} }.
$
One can easily verify that
$$
\dfrac {|y-l|^2+ 1 } {|y|^2+1}\leq (|l|+1)^2,\quad \text{for $k>0$}
$$
and hence \fer{sepx} follows.\\

\noindent The above bound on $s_{\e,x}^i$ helps us to control $g$. We obtain indeed, for some positive constant $C$,
$$
\begin{array}{rl}
g&= \displaystyle\int_{k\geq 1} \left[s_{\e,x}^1(k)
+ s_{\e,x}^2(k)  -2\right] \f{e^k}{|e^k-1|^{1+\al}}dk 
\,\leq \,
2\displaystyle\int_{k\geq 1} e^{Ak}\f{e^k}{|e^k-1|^{1+\al}}dk \leq C.
\end{array}
$$
Note that the above integral is bounded since $A<\al$.\\

\noindent To control $f$, we compute the Taylor expansion of $s_{\e,x}^1+s_{\e,x}^2$ around $k=0$:
$$
s_{\e,x}^1(k)+s_{\e,x}^2(k)=2+\f{k^2}{2}\f{d^2}{dk^2}\left( s_{\e,x}^1+s_{\e,x}^2\right)(k'),\quad \text{with $0\leq k'\leq k\leq 1$.}
$$
Using \fer{sepx}, it is easy to show that for $0\leq \e \leq \e_0$ and $0\leq k'\leq 1$, we have
$$
| \f{d^2}{dk^2}\left( s_{\e,x}^1+s_{\e,x}^2\right)(k') |\leq C_0(\e_0),
$$
where $C_0$ is a positive constant depending only on $\e_0$. It then follows  that, for a large enough constant $C$,
$$
| f |\leq \f{C(\e_0)}{2}\int_{0\leq k\leq 1} k^2 \f{e^k}{|e^k-1|^{1+\al}}dk \leq C,
$$
since  $\al<2$ and $e^k-1\approx k$ near $k=0$.\\

\noindent Combining the above bounds on $f$ and $g$, we obtain that for  large enough constant $C$ and for all $\e\leq \e_0$,
$\
S\leq C.
$
\\

(ii)  \textbf{Uniform bounds from below.}  We prove that for $D$ and $C$ large enough constants,
$$
 \underline s(x,t):=-\f A 2\log(|x|^2+1)-D-Ct\leq u_\e(x,t),\qquad \text{for all $(x,t)\in \R\times \R^+$}.
$$
We first prove that the above inequality is verified for $t=0$, for $D$ large enough. We then show that, for $C$ large enough,  $\underline s$ is a subsolution of \fer{eq:ue},   where we fix the last term $R(I_\e)$, with $I_\e=\int e^{\f{u_\e}{\e}}dx$. Then, the claim follows from the comparison principle since $u_\e$ is a solution and in particular a supersolution of \fer{eq:ue}.

\noindent To prove the inequality for $t=0$, we first notice from \fer{as:biz} that
$$
u_\e^0(0)-A\log \left( 1+ |x| \right) \leq u_\e^0(x),\qquad \text{for all $x\in \R$.}
$$
Next, we notice from \fer{as:ue0} that for $\e_0$ small enough, $u_\e^0(0)$ is uniformly bounded for $0\leq \e\leq \e_0$. 
Therefore, we can choose $D$ large enough, such that for $\e\leq \e_0$,
$$
-D- \f A2 \log \left( 1+ |x|^2 \right) \leq u_\e^0(0)-A\log \left( 1+ |x| \right) \leq u_\e^0(x),\qquad \text{for all $x\in \R$.}
$$

\noindent To prove that $\underline s$ is a subsolution of \fer{eq:ue},  since $R$ is bounded thanks to \fer{as:monR} and \fer{boundI}, it is enough to prove that, for $C$ large enough,
 $$
 S=\int_{k\geq 0} \left[ e^{\f{\overline s(x+e^{\e k}-1,t)-\overline s(x,t)}{\e}}+e^{\f{\overline s(x-e^{\e k}+1,t)-\overline s(x,t)}{\e}}-2\right] \f{e^k}{|e^k-1|^{1+\al}}dk\geq -C.
$$
As in Step (i) above we split $S$ into two terms $S=f+g$, with $f$ and $g$ given respectively by \fer{def:f} and \fer{def:g}.
The term $f$ can be controlled in the same way as in Step (i) in the proof of Theorem \ref{th:reg}. To control $g$ we compute 
$$
\displaystyle\int_{k\geq 1} \left[ \left( \dfrac { |x|^2+1 }  { |x+e^{\e k}-1|^2+1 }  \right) ^{\f{A}{2\e} } 
+  \left( \dfrac { |x|^2+1 }  { |x-e^{\e k}+1|^2+1 }  \right) ^{\f{A}{2\e} }  -2\right] \f{e^k}{|e^k-1|^{1+\al}}dk\geq
-2\displaystyle\int_{k\geq 1} \f{e^k}{|e^k-1|^{1+\al}}dk,
$$
which is enough to conclude, since the r.h.s. of the above inequality is bounded from below.
\\

(iii) \textbf{The proof of \fer{log-es}.}\\

For all $h\in \R$ and $\e>0$, we define 
$$
w_{\e,h}(x,t)=u_\e(x+h,t)-u_\e(x,t),\quad \text{for $t\geq 0$ and $x\in \R$.}
$$
We then compute 
$$
\begin{array}{rl}
\p_t w_{\e,h}(x,t) &=\int_{k\geq 0} \left[ e^{\f{u_\e(x+h+e^{\e k}-1,t)-u_\e(x+h,t)}{\e}} - e^{\f{u_\e(x+e^{\e k}-1,t)-u_\e(x,t)}{\e}} \right. \\
&\left. + e^{\f{u_\e(x+h - e^{\e k} +1,t)-u_\e(x+h,t)}{\e}}- e^{\f{u_\e(x-e^{\e k}+1,t)-u_\e(x,t)}{\e}}  \right]  \f{e^k}{|e^k-1|^{1+\al}}dk.
\end{array}
$$
Using the convexity inequality
$\
e^a\leq e^b+e^a(a-b),
$
we deduce that
$$
\begin{array}{rl}
\p_t w_{\e,h}(x,t) &\leq \int_{k\geq 0} \left[ e^{\f{u_\e(x+h+e^{\e k}-1,t)-u_\e(x+h,t)}{\e}} \left(\f{w_{\e,h}(x+e^{\e k}-1,t)-w_{\e,h}(x,t)}{\e}\right) \right. \\
&\left. + e^{\f{u_\e(x+h - e^{\e k} +1,t)-u_\e(x+h,t)}{\e}} \left(\f{w_{\e,h}(x-e^{\e k}+1,t)-w_{\e,h}(x,t)}{\e}\right)  \right]  \f{e^k}{|e^k-1|^{1+\al}}dk.
\end{array}
$$
Therefore, by the maximum principle, \fer{boundue} and  \fer{as:biz} we obtain that for all $t>0$, $\e>0$ and $h,\,x\in\R$,
$$
w_{\e,h}(x,t)\leq {\sup_x\, }w_{\e,h}(x,0) \leq A\log (1+|h|),
$$
and hence \fer{log-es} follows.\\

(iv) \textbf{The proof of \fer{maxue}.} We prove \fer{maxue}, we first notice from \fer{boundI} that
$$
0<I_m\leq \int_{\R} e^{\f{u_\e(x,t)}{\e}}dx \leq I_M.
$$
Moreover, we already know from step (i) that
$\
u_\e(x,t) \leq-\f A 2\log(|x|^2+1)+B+Ct.
$

\noindent The two above properties imply that there exists $A_2=A_2(T)$ large enough such that, for all $t\in [0,T]$ and $\e \leq \e_0$ with $\e_0=\e_0(A)$ small enough,
$$
\f {I_m}{2}\leq \int_{|x| \leq A_2} e^{\f{u_\e(x,t)}{\e}}dx.
$$
We deduce that 
$
\  \e \log \left( \f{I_m}{4A_2(T)} \right) \leq \max_{x\in B_{A_2}(0)} u_\e(x,t),\  \text{for all $t\in  [0,T]$ and $\e\leq \e_0$},
$
and hence \fer{maxue}. \\
\qed
\section{Regularity results for \fer{SMe} and the reaction term given by \fer{re-2}}
\label{sec:reg2}

\begin{theorem}
\label{th:reg2}
Assume  \fer{as:ue0}, \fer{as:u0}, \fer{as:biz}  and \fer{as:n01}. Then, for all $T>0$ and $R>0$, there exist constants $\e_0$, $A_1(R,T)$, $A_2(T)$, $D$ and $C$ such that, for all $\e \leq \e_0$,
\beq
\label{boundue2}
 -\f A 2\log(|x|^2+1)-D-Ct \leq u_\e\leq -\f A 2\log(|x|^2+1)+B+Ct,\quad \text{in $ B_R(0)\times [0,T]$},
 \eeq
Moreover, we have
\beq
\label{log-es2}
u_\e(x+h,t) \leq u_\e(x,t)+A\log(1+|h|),\qquad \text{for all $x,\, h \in \R$ and $t\geq 0$}.
\eeq
In particular $(u_\e)_\e$ is uniformly Lipschitz with respect to $x$:
\beq
\label{Lipu2}
\| D_x u_\e \|_{L^\infty(\R\times \R^+)}\leq A.
\eeq


\end{theorem}

\proof
(i) \textbf{Uniform bounds from above and below.} The inequalities given in \fer{boundue2} can be proved following similar arguments as in the proof of Steps (i) and (ii) in Theorem \ref{th:reg}. The only difference here is that the boundedness of the reaction term $R$ is derived from \fer{boundn}.\\

(ii) \textbf{The proof of \fer{log-es2}.} The proof of this part is also close to the one in Theorem \ref{th:reg}.  As in the previous case, for all $h\in \R$ and $\e>0$, we define 
$$
w_{\e,h}(x,t)=u_\e(x+h,t)-u_\e(x,t),\quad \text{for $t\geq 0$ and $x\in \R$.}
$$
We then compute 
$$
\begin{array}{rl}
\p_t w_{\e,h}(x,t) &=\int_{k\geq 0} \left[ e^{\f{u_\e(x+h+e^{\e k}-1)-u_\e(x+h)}{\e}} - e^{\f{u_\e(x+e^{\e k}-1)-u_\e(x)}{\e}} \right. \\
&\left. + e^{\f{u_\e(x+h - e^{\e k} +1)-u_\e(x+h)}{\e}}- e^{\f{u_\e(x-e^{\e k}+1)-u_\e(x)}{\e}}  \right] \f{e^k}{|e^k-1|^{1+\al}} dk+n_\e(x,t)-n_\e(x+h,t).
\end{array}
$$
Using the convexity inequality
$\
e^a\leq e^b+e^a(a-b),
$
we deduce that
$$
\begin{array}{rl}
\p_t w_{\e,h}(x,t) &\leq \int_{k\geq 0} \left[ e^{\f{u_\e(x+h+e^{\e k}-1)-u_\e(x+h)}{\e}} \left(\f{w_{\e,h}(x+e^{\e k}-1)-w_{\e,h}(x)}{\e}\right) \right. \\
&\left. + e^{\f{u_\e(x+h - e^{\e k} +1)-u_\e(x+h)}{\e}} \left(\f{w_{\e,h}(x-e^{\e k}+1)-w_{\e,h}(x)}{\e}\right)  \right] \f{e^k}{|e^k-1|^{1+\al}} dk+n_\e(x,t)-n_\e(x+h,t).
\end{array}
$$
Therefore, by the maximum principle, \fer{as:biz}, \fer{boundue2} and since $u_\e(x+h,t)-u_\e(x,t)$ and $n_\e(x+h,t)-n_\e(x,t)$ have the same sign, we obtain that for all $t>0$, $\e>0$ and $h,\,x\in\R$,
$$
w_{\e,h}(x,t)\leq \max \left(0, \sup_x w_{\e,h}(x,0) \right) \leq A\log (1+|h|).
$$
and hence \fer{log-es2} follows.\\
\qed

\section{Proof of Theorem \ref{th:main}}
\label{sec:th1}
To prove Theorem \ref{th:main}, we use the half-relaxed methods for viscosity solutions \cite{C.I.L:92,GB:94}. 
Since $(u_\e)_\e$ is locally uniformly bounded, we can define it's lower and upper semicontinuous envelopes
\begin{equation*} 
\underline u(x,t) := \underset{\underset{(y,s)\to (x,t)}{\e\to 0}}{ \underline \liminf} u_\e(y,s),
\qquad
\overline u(x,t) := \underset{\underset{(y,s)\to (x,t)}{\e\to 0}}{ \overline \limsup} u_\e(y,s).
\end{equation*}

(i) We prove Theorem \ref{th:main}-(i), in several steps. We first prove that $\underline u$ is a viscosity supersolution of \fer{HJ}. Next we prove that $\overline u$ is a viscosity subsolution of \fer{HJ}. We then conclude using that \fer{HJ} admits a comparison principle. Finally we prove \fer{max}. \\

{\bf Step 1.} {\bf ($\underline u$ is a viscosity supersolution of \fer{HJ})} Let $\vp\in \mathcal{C} \left( \R \times \R^+ \right) \cap \mathcal{C}^2  \left(\Omega(x_0,t_0)\right)$, with $\Omega(t_0,x_0)$ an open neighborhood of $(x_0,t_0)$, be a test function. We assume that $\underline u-\vp$ has a global minimum at $(x_0,t_0)$. By classical arguments in the theory of viscosity solutions (see \cite{C.I.L:92,GB:94}) we can assume that the minimum at $(x_0,t_0)$ is strict and thus there exists a sequence $(x_\e,t_\e)$ such that $(x_\e,t_\e)$ tends to $(x_0,t_0)$, and  $u_\e(x_\e,t_\e)$ tends to $\underline u(x_0,t_0)$ as $\e\to 0$ and $u_\e-\vp$ takes a minimum at $(x_\e,t_\e)$. Since $u_\e$ solves \fer{eq:ue}, we find
$$
\p_t \vp(x_\e,t_\e) - R(I_\e(t_\e))\geq \int_{k\geq 0} \left[ e^{\f{u_\e(x_\e+e^{\e k}-1,t_\e)-u_\e(x_\e,t_\e)}{\e}}+e^{\f{u_\e(x_\e-e^{\e k}+1,t_\e)-u_\e(x_\e,t_\e)}{\e}}-2\right] \f{e^k}{|e^k-1|^{1+\al}}dk.
$$
Since $u_\e-\vp$ takes a minimum at $(x_\e,t_\e)$, we obtain
$$
\vp(x_\e+l,t_\e)-\vp(x_\e,t_\e)\leq u(x_\e+l,t_\e) - u(x_\e,t_\e),\qquad \text{for all $l\in \R$}.
$$
It follows that
\beq
\label{eq:int}
\begin{array}{rl}
\p_t \vp(x_\e,t_\e)\geq &  R(I_\e(t_\e))
\\
+&\int_{M\geq k\geq 0} \left[e^{\f{\vp(x_\e+e^{\e k}-1,t_\e)-\vp(x_\e,t_\e)}{\e}}+e^{\f{\vp(x_\e-e^{\e k}+1,t_\e)-\vp(x_\e,t_\e)}{\e}}-2\right] \f{e^k}{|e^k-1|^{1+\al}}dk\\
+&\int_{k\geq M} \left[ e^{\f{u_\e(x_\e+e^{\e k}-1,t_\e)-u_\e(x_\e,t_\e)}{\e}}+e^{\f{u_\e(x_\e-e^{\e k}+1,t_\e)-u_\e(x_\e,t_\e)}{\e}}-2\right]  \f{e^k}{|e^k-1|^{1+\al}}dk.
\end{array}
\eeq
We note that, using the Taylor-Lagrange formula, for some $0< \mu,\,\mu' < \e$
$$
\begin{array}{rl}
\f{\vp(x_\e+e^{\e k}-1,t_\e)-\vp(x_\e,t_\e)}{\e}&=D_x \vp(x_\e,t_\e)\cdot k\\
& +\f{\e}{2} \left[e^{\mu k}k^2 D_x \vp(x_\e+e^{\mu k}-1,t_\e)+e^{2\mu k}k^2 D^2 \vp(x_\e+e^{\mu k}-1,t_\e)\right],
\end{array}
$$
$$
\begin{array}{rl}
\f{\vp(x_\e-e^{\e k}+1,t_\e)-\vp(x_\e,t_\e)}{\e}&=-D_x \vp(x_\e,t_\e)\cdot k\\
& +\f{\e}{2} \left[-e^{\mu' k}k^2 D_x \vp(x_\e-e^{\mu' k}+1,t_\e)+e^{2\mu' k}k^2 D^2 \vp(x_\e-e^{\mu' k}+1,t_\e)\right].
\end{array}
$$
Since $\vp\in \mathcal{C}^2  \left(\Omega(x_0,t_0)\right)$, it follows that, for fixed $M$ and as $\e\to 0$, the second term of the r.h.s. of \fer{eq:int} converges to
$$
\int_{M\geq k\geq 0} \left[e^{D_x \vp(x_0,t_0)\cdot k}+e^{-D_x \vp(x_0,t_0)\cdot k}-2\right] \f{e^k}{|e^k-1|^{1+\al}}dk.
$$
Furthermore, one can control the third term of the r.h.s. of \fer{eq:int} as below
$$
\begin{array}{c}
\displaystyle\int_{k\geq M} \left[ e^{\f{u_\e(x_\e+e^{\e k}-1,t_\e)-u_\e(x_\e,t_\e)}{\e}}+e^{\f{u_\e(x_\e-e^{\e k}+1,t_\e)-u_\e(x_\e,t_\e)}{\e}}-2\right]  \f{e^k}{|e^k-1|^{1+\al}}dk\
\geq \ -2\int_{k\geq M} \f{e^k}{|e^k-1|^{1+\al}}dk.
\end{array}
$$
Combining the above lines and Lemma \ref{th:I} we deduce
$$
\begin{array}{c}
\p_t \vp(x_0,t_0)\geq\,
\displaystyle\int_{M\geq k\geq 0} \left[e^{D_x \vp(x_0,t_0)\cdot k}+e^{-D_x \vp(x_0,t_0)\cdot k}-2\right] \f{e^k}{|e^k-1|^{1+\al}}dk\
-2 \displaystyle\int_{k\geq M} \f{e^k}{|e^k-1|^{1+\al}}dk.
\end{array}
$$
Letting $M\to \infty$ we obtain
$$
\p_t \vp(x_0,t_0)\geq
\int_{ k\geq 0} \left[e^{D_x \vp(x_0,t_0)\cdot k}+e^{-D_x \vp(x_0,t_0)\cdot k}-2\right] \f{e^k}{|e^k-1|^{1+\al}}dk.
$$
It follows that $\underline u$ is a viscosity supersolution of \fer{HJ}.\\

{\bf Step 2.} {\bf ($\overline u$ is a viscosity subsolution of \fer{HJ})} Let $\vp\in \mathcal{C} \left( \R \times \R^+ \right) \cap \mathcal{C}^2  \left(\Omega(x_0,t_0)\right)$, with $\Omega(t_0,x_0)$ an open neighborhood of $(x_0,t_0)$,  be a test function. We assume that $\overline u-\vp$ has a global maximum at $(x_0,t_0)$. We prove that 
\beq
\label{eq:sub}
\p_t \vp(x_0,t_0)\leq
\int_{ k\geq 0} \left[e^{D_x \vp(x_0,t_0)\cdot k}+e^{-D_x \vp(x_0,t_0)\cdot k}-2\right] \f{e^k}{|e^k-1|^{1+\al}}dk.
\eeq
We first notice from \fer{Lipu} that 
$$
|D_x \vp|(x_0,t_0) \leq A<\al.
$$
By similar arguments as in the previous steps, we obtain that there exist a sequence $(x_\e,t_\e)$ such that $u_\e-\vp$ takes a maximum at $(x_\e,t_\e)$ and that
$$
\begin{array}{rl}
\p_t \vp(x_\e,t_\e)\leq &  R(I_\e(t_\e))
\\
+&\int_{ M\geq k\geq 0} \left[e^{\f{\vp(x_\e+e^{\e k}-1,t_\e)-\vp(x_\e,t_\e)}{\e}}+e^{\f{\vp(x_\e-e^{\e k}+1,t_\e)-\vp(x_\e,t_\e)}{\e}}-2\right] \f{e^k}{|e^k-1|^{1+\al}}dk\\
+&\int_{  k\geq M} \left[e^{\f{u_\e(x_\e+e^{\e k}-1,t_\e)-u_\e(x_\e,t_\e)}{\e}}+e^{\f{u_\e(x_\e-e^{\e k}+1,t_\e)-u_\e(x_\e,t_\e)}{\e}}-2\right] \f{e^k}{|e^k-1|^{1+\al}}dk.
\end{array}
$$
Again following similar arguments as above, the second term of the r.h.s. of the above inequality converges to
$$
\int_{M\geq k\geq 0} \left[e^{D_x \vp(x_0,t_0)\cdot k}+e^{-D_x \vp(x_0,t_0)\cdot k}-2\right] \f{e^k}{|e^k-1|^{1+\al}}dk.
$$
Moreover, from \fer{log-es} we obtain
$$
\begin{array}{rl}
&\displaystyle\int_{  k\geq M} \left[e^{\f{u_\e(x_\e+e^{\e k}-1,t_\e)-u_\e(x_\e,t_\e)}{\e}}+e^{\f{u_\e(x_\e-e^{\e k}+1,t_\e)-u_\e(x_\e,t_\e)}{\e}}-2\right] \f{e^k}{|e^k-1|^{1+\al}}dk\\
&\leq  \displaystyle\int_{  k\geq M} \left[2 e^{\f{A\log(1+e^{\e k}-1)}{\e}}-2\right] \f{e^k}{|e^k-1|^{1+\al}}dk\
\leq\  2 \displaystyle\int_{  k\geq M} \f{e^{(A+1)k}}{|e^k-1|^{1+\al}}dk.
\end{array}
$$
Combining the above arguments and Theorem \ref{th:I} we deduce that
$$
\begin{array}{rl}
\p_t \vp(x_0,t_0)\leq &  
\displaystyle\int_{M\geq k\geq 0} \left[e^{D_x \vp(x_0,t_0)\cdot k}+e^{-D_x \vp(x_0,t_0)\cdot k}-2\right] \f{e^k}{|e^k-1|^{1+\al}}dk\
+2\displaystyle\int_{  k\geq M} \f{e^{(A+1)k}}{|e^k-1|^{1+\al}}dk.
\end{array}
$$
Letting $M$ go to infinity, and in view of $A<\al$, we obtain
$$
\p_t \vp(x_0,t_0)\leq
\int_{ k\geq 0} \left[e^{D_x \vp(x_0,t_0)\cdot k}+e^{-D_x \vp(x_0,t_0)\cdot k}-2\right] \f{e^k}{|e^k-1|^{1+\al}}dk.
$$

{\bf Step 3.} {\bf (Convergence of  $(u_\e)_\e$ to the unique solution of \fer{HJ})} From the above steps we obtain that $\underline u$  and $\overline u$ are respectively  viscosity supersolution and  viscosity subsolution of \fer{HJ}. Moreover, combing the above arguments with \fer{as:ue0}, we also obtain that $\underline u$  and $\overline u$ are  viscosity supersolution and  viscosity subsolution of  \fer{HJ} up to the boundary $\R \times \{0\}$.   Finally, in the one hand, from the strong comparison principle satisfied by \fer{HJ} (see for instance \cite{GB:94}), we obtain that 
$\
\overline u \leq \underline u.
$
In the other hand, by definition we also have
$\ 
\underline u \leq \overline u.
$
It follows that $(u_\e)_\e$ converges locally uniformly to $u=\underline u=\overline u$. \\

{\bf Step 4.} {\bf (Proof of \fer{max})} Firstly, the first part of \fer{max} is a consequence of \fer{Lipu} and the uniform convergence of $(u_\e)_\e$ to $u$. We next deduce from \fer{maxue} that 
$\
0 \leq \max_{x\in \R} u(x,t), \ \text{for all $t\in \R^+$}.
$
Finally, we obtain from the upper bound in \fer{boundI} and  the first part of \fer{max} that 
$\
 \max_{x\in \R} u(x,t)\leq 0, \ \text{for all $t\in \R^+$},
$
and hence the second part of \fer{max}. \\

(ii) We first deduce from \fer{boundI} that, along subsequences as $\e\to 0$, $n_\e$ converges in $L^\infty\left( w*(0,\infty) ; \mathcal{M}^1(\R) \right)$ to a measure $n$. Next, we use \fer{Hopf} and the fact that $(u_\e)_\e$ converges locally uniformly to $u$ to obtain that,  $supp\, n \subset \{(x,t)\,|\, u(x,t)=0\}$. 
\section{Proof of Theorem \ref{th:main2}}
\label{sec:th2}
To prove Theorem \ref{th:main2}, we use the same scheme as in Section \ref{sec:th1}. We first prove that $\underline u$ is a viscosity  supersolution of \fer{HJ2}. Next we prove that $\overline u$ is a viscosity subsolution of \fer{HJ2}. Next, noticing that \fer{HJ2} admits a comparison principle (see for instance \cite{GB:94}  and \cite{LE.PS:89}), we conclude that $(u_\e)_\e$ converges locally uniformly to the unique viscosity solution of \fer{HJ2}. Furthermore,   \fer{LipKPP} is a consequence of  \fer{Lipu2} and the uniform convergence of $(u_\e)_\e$ to $u$. Finally we prove \fer{limn2}.\\

{\bf Step 1.} {\bf ($\underline u$ is a viscosity supersolution of \fer{HJ2})} 
We first notice that  if $\underline  u(x_0,t_0)\geq 0$, the supersolution criterion for \fer{HJ2} is obviously verified at $(x_0,t_0)$. Therefore,  it is enough to study only the case $\underline u(x_0,t_0)<0$.

 \noindent Let $\vp\in \mathcal{C} \left( \R \times \R^+ \right) \cap \mathcal{C}^2  \left(\Omega(x_0,t_0)\right)$, with $\Omega(t_0,x_0)$ an open neighborhood of $(x_0,t_0)$, be a test function. We assume that $\underline u-\vp$ has a global minimum at $(x_0,t_0)$. As previously, by classical arguments in the theory of viscosity solutions  we can assume that the minimum at $(x_0,t_0)$ is strict and thus there exist a sequence $(x_\e,t_\e)$ such that $(x_\e,t_\e)$ tends to $(x_0,t_0)$, and  $u_\e(x_\e,t_\e)$ tends to $\underline u(x_0,t_0)$ as $\e\to 0$ and $u_\e-\vp$ takes a minimum at $(x_\e,t_\e)$. Since $u_\e$ solves \fer{eq:ue}, we find
$$
\p_t \vp(x_\e,t_\e) -1+ e^{\f{u_\e(x_\e,t_\e)}{\e}}\geq \int_{k\geq 0} \left[ e^{\f{u_\e(x_\e+e^{\e k}-1,t_\e)-u_\e(x_\e,t_\e)}{\e}}+e^{\f{u_\e(x_\e-e^{\e k}+1,t_\e)-u_\e(x_\e,t_\e)}{\e}}-2\right] \f{e^k}{|e^k-1|^{1+\al}}dk.
$$
We then deduce, following similar arguments as in Step (i) in Section \ref{sec:th1}, that
$$
\p_t \vp(x_0,t_0)-
\int_{ k\geq 0} \left[e^{D_x \vp(x_0,t_0)\cdot k}+e^{-D_x \vp(x_0,t_0)\cdot k}-2\right] \f{e^k}{|e^k-1|^{1+\al}}dk \geq \limsup_{\e\to 0} \left(1- e^{\f{u_\e(x_\e,t_\e)}{\e}} \right).
$$
Moreover, since $u_\e(x_\e,t_\e)$ tends to $\underline u(x_0,t_0)$ as $\e\to 0$ and $\underline  u(x_0,t_0)<0$, the r.h.s. of the above inequality is equal to $1$. We deduce that
$$
\p_t \vp(x_0,t_0)-
\int_{ k\geq 0} \left[e^{D_x \vp(x_0,t_0)\cdot k}+e^{-D_x \vp(x_0,t_0)\cdot k}-2\right] \f{e^k}{|e^k-1|^{1+\al}}dk -1 \geq 0.
$$

{\bf Step 2.} {\bf ($\overline u$ is a viscosity subsolution of \fer{HJ2})}
We first notice from \fer{boundn} that $\overline u(x,t)\leq 0$, for all $(x,t)\in \R\times \R^+$. Therefore, it is enough to prove that $\overline u$ is a viscosity subsolution of 
$$
\p_t u-\int_0^\infty \left( e^{D_x  u\cdot k}+e^{-D_x  u\cdot k}-2 \right) \f{e^kdk}{|e^k-1|^{1+\al}}-1\leq 0.
$$\\
Let $\vp\in \mathcal{C} \left( \R \times \R^+ \right) \cap \mathcal{C}^2  \left(\Omega(x_0,t_0)\right)$, with $\Omega(t_0,x_0)$ an open neighborhood of $(x_0,t_0)$,  be a test function. We assume that $\overline u-\vp$ has a global maximum at $(x_0,t_0)$, which implies as previously that there exist a sequence $(x_\e,t_\e)$ such that $(x_\e,t_\e)$ tends to $(x_0,t_0)$ and  $u_\e(x_\e,t_\e)$ tends to $\underline u(x_0,t_0)$ as $\e\to 0$, and $u_\e-\vp$ takes a maximum at $(x_\e,t_\e)$. We deduce that
$$
\begin{array}{rl}
\p_t \vp(x_\e,t_\e)\leq &  1-n_\e(x_\e,t_\e)+\int_{ k\geq 0} \left[e^{\f{\vp(x_\e+e^{\e k}-1,t_\e)-\vp(x_\e,t_\e)}{\e}}+e^{\f{\vp(x_\e-e^{\e k}+1,t_\e)-\vp(x_\e,t_\e)}{\e}}-2\right] \f{e^k}{|e^k-1|^{1+\al}}dk.
\\
\leq&1+\int_{ k\geq 0} \left[e^{\f{\vp(x_\e+e^{\e k}-1,t_\e)-\vp(x_\e,t_\e)}{\e}}+e^{\f{\vp(x_\e-e^{\e k}+1,t_\e)-\vp(x_\e,t_\e)}{\e}}-2\right] \f{e^k}{|e^k-1|^{1+\al}}dk.
\end{array}
$$
It then follows following similar arguments as in Step (i) in Section \ref{sec:th1}, that
$$
\p_t \vp(x_0,t_0)\leq
1+
\int_{ k\geq 0} \left[e^{D_x \vp(x_0,t_0)\cdot k}+e^{-D_x \vp(x_0,t_0)\cdot k}-2\right] \f{e^k}{|e^k-1|^{1+\al}}dk,
$$
and hence $\overline u$ is a viscosity subsolution of \fer{HJ2}.\\

{\bf Step 3.} {\bf (The proof of \fer{limn2})} Let $(x_0,t_0)$ be such that $u(x_0,t_0)<0$. It follows easily from \fer{Hopf} and the locally uniform convergence of $(u_\e)$ to $u$, that $n_\e$ goes to $0$ locally uniformly, as $\e\to 0$.

\noindent We now suppose that, there exists $r,\, \da>0$ such that
$(\tilde x-2r,\tilde x+2r)\times (\tilde  t-2\da,\tilde t+2\da) \subset  \{(x,t)\in \R\times (0,\infty) \, |\,u(x,t)=0 \}$. Let $(x_0,t_0)\in (\tilde x-r,\tilde x+r)\times (\tilde  t-\da,\tilde t+\da)$.  We consider the following test function:
$$
\vp(x,t) = -\f A r (x-x_0)^2-(t-t_0)^2.
$$
One can verify easily that $u-\vp$ has a local minimum at $(x_0,t_0)$. We show that  this minimum point is indeed global with respect to $x$. We first find from \fer{log-es2} that
$$
-A \log \left( 1+ |x-x_0| \right) \leq u(x,t), \qquad \text{for all $(x,t)\in \R\times  (t_0-\da,t_0+\da)$}.
$$
Next, we notice that
$$
 -\f A r (x-x_0)^2 < -A \log \left( 1+ |x-x_0| \right), \qquad \text{for all $|x-x_0|>r$}.
$$
Combining the above inequalities and the fact that $(x_0-r,x_0+r)\times (t_0-\da,t_0+\da) \subset  \{(x,t)\in \R\times (0,\infty) \, |\,u(t,x)=0 \}$, we deduce that $u-\vp$ has a  minimum at $(x_0,t_0)$ which is global with respect to $x$. Moreover, this is a strict minimum. It follows that  there exist points $(x_\e,t_\e)\in (x_0-r,x_0+r)\times (t_0-\da,t_0+\da)$ such that  $u_\e-\vp_\e$ has a  local in $t$ and global in $x$ minimum at $(x_\e,t_\e)$ and such that $(x_\e,t_\e)\to (x_0,t_0)$.\\

Since $u_\e -\vp$ has a  local in $t$ and global in $x$ minimum at $(x_\e, t_\e)$, we have
 $$
  \p_t u_\e (x_\e,t_\e)= \p_t \vp_\e (x_\e,t_\e)=-2(t_\e-t_0),
 $$
$$
\begin{array}{c}
\displaystyle\int_0^\infty \left( e^{\f{u_\e \left( x_\e+e^{\e k} -1,t_\e \right)- u_\e(x_\e,t_\e)}{\e} }+ e^{\f{u_\e \left(x_\e-e^{\e k} +1,t_\e \right)- u_\e(x_\e,t_\e)}{\e} }  -2 \right) \f{e^kdk}{|e^k-1|^{1+\al}}\\
\geq
\displaystyle\int_0^\infty
 \left( e^{\f{A\left( (x_\e-x_0)^2 - ( x_\e-x_0+e^{\e k} -1)^2 \right)}{r\e}} + e^{\f{A\left( (x_\e-x_0)^2 - ( x_\e-x_0-e^{\e k} +1)^2 \right)}{r\e} }  -2 \right)
  \f{e^kdk}{|e^k-1|^{1+\al}}\geq o(1).
 \end{array}
 $$
Combining the above lines with \fer{eq:ue} we deduce that
$\
n_\e(x_\e,t_\e)\geq 1+o(1).
$
Moreover, following similar arguments as in the proof of Theorem \ref{th:kpp}, part (ii), we obtain that
$\ 
n_\e(x_0,t_0) \geq n_\e(x_\e,t_\e),
$
and hence 
$$
\liminf_{\e\to 0} n_\e(x_0,t_0) \geq 1,\qquad \text{uniformly in $(x_0-r,x_0+r)\times (t_0-\da,t_0+\da)$}.
$$
Finally, we conclude from the above inequality and Lemma \ref{th:n} that $n_\e(x_0,t_0)\to 1$ uniformly in $(x_0-r,x_0+r)\times (t_0-\da,t_0+\da)$, as $\e\to 0$.

\section{The multi-dimensional case}
\label{sec:multid}

In this section we show how the above results can be generalized to the multidimensional case $x\in \R^N$.

\subsection{The long range/long time rescaling}


To introduce the rescaling for the multidimensional case, we define the following mapping
$$
p(z)=
\begin{cases}
\f{z}{|z|}& \text{for $z\in \R^N \setminus \{0\}$,}\\
0& \text{for $z=0$.}
\end{cases}
$$
We then introduce the following rescaling
$$
x\mapsto |x|^{\f{1}{\e}} p(x), \qquad t\mapsto \f{t}{\e},
\qquad
n_\e(x,t) =n\big(|x|^{\f 1 \e}\, p(x), \f t \e\big).
$$
We replace this in \fer{SM} with $x\in \R^N$, and obtain, 
\beq
\label{KPPe-d}
\begin{cases}
\e \p_t n_\e (x,t) =\int_0^\infty \int_{\nu\in S^{N-1}} \left( n_\e \left( \left| |x|^{\f 1\e} p(x)+ h\nu \right|^\e \,p( |x|^{\f 1\e}p(x)+h\nu),t \right)  -n_\e(x,t) \right) \f{dS\,dh}{|h|^{1+\al}}\\
\hskip 2cm + \,n_\e(x,t) R(n_\e,I_\e)(x,t),\\
n_e(x,0)=n_\e^0(x),
\end{cases}
\eeq
where $I_\e(t)=I(\f t \e)$. With this rescaling, we can obtain the macroscopic behavior of the dynamics as before and extend Theorems \ref{th:kpp} and \ref{th:RI}  to the case with $x\in \R^N$:

\begin{theorem}\label{th:kpp-N}
Let $x\in \R^N$ and  $n_\e$ be the solution of \fer{KPPe-d} with \fer{re-2} and $u_\e=\e\log n_\e$. 

(i) Under assumption \fer{as:inikpp}, as $\e\to 0$, $(u_\e)_\e$ converges locally uniformly to $u$ defined as below
$$
u(x,t)=\min( 0, -(1+\al) \log |x|+t).
$$
(ii) Moreover, as $\e\to 0$,
$$
\begin{cases}
n_\e\to 0, & \text{locally uniformly in  $\mathcal A=\{(x,t)\in \R^N\times (0,\infty) \, |\, t <(1+\al)\log |x|\}$,}\\
 n_\e \to 1, & \text{locally uniformly in  $\mathcal B=\{(x,t)\in \R^N\times (0,\infty) \, |\, t >(1+\al)\log |x|\}$.}
\end{cases}
$$
\end{theorem}

\begin{theorem}
\label{th:RI-N}
Let $x\in \R^N$ and $n_\e$ be the solution of \fer{KPPe-d} with \fer{re-1} and $u_\e=\e\log n_\e$. 

(i) Under assumptions  \fer{as:inikpp},  \fer{as:maxR}, \fer{as:monR} and \fer{as:I}, as $\e\to 0$, $(u_\e)_\e$ converges locally uniformly to $u\in \mathcal{C}(\R^N)$ defined as below
$$
u(x,t)=\min(0,-(1+\al) \log |x|).
$$
(ii) Moreover, $n_\e$ converges, along subsequences as $\e\to 0$, in $L^\infty$ weak-$\ast$ to a function $n\in L^\infty(\R^N\times \R^+)$, such that $\mathrm{supp} \; n\subset \{(x,t)\in \R^N\times \R^+ \, | \, u(x,t)=0\}= \{(x,t)\in  \R^N\times \R^+ \,| \, |x|\leq 1\} $.
\end{theorem}

\proof[Proof of Theorems \ref{th:kpp} and \ref{th:RI}]
Note that the proofs of Theorems \ref{th:kpp} and \ref{th:RI} are based on Lemma \ref{lem:ex}. We claim that an equivalent lemma holds in the multidimensional case. 

\begin{lemma}
\label{lem:exN}
Let $g_N:\R^N\to \R$ be given by $g_N(x)=\f{1}{1+|x|^{1+\al}}$. Then, there exists a positive constant $C_N$, independent of $x$, such that
\beq
\label{gN}
|(-\Delta)_N^{\f\al 2} g_N(x)|\leq C_N g_N(x),
\eeq
where $ (-\Delta )_N^{\f \al 2}$ is the $N$-dimensional fractional laplacian, such that
\beq
\label{int-N}
\f{(-\Delta )_N^{\f \al 2}g_N(x)}{g_N(x)}  =\int_0^\infty \int_{\nu\in S^{N-1}} \left( \f{1+|x|^{1+\al}}{1+|x+h\nu|^{1+\al}}-1\right)\f{dSdh}{|h|^{1+\al}}.
\eeq
\end{lemma}

\noindent One can easily verify that, replacing the result of Lemma \ref{lem:ex} by Lemma \ref{lem:exN},  the other parts of the proofs will be easily adapted for $x\in \R^N$. We prove Lemma \ref{lem:exN} in Appendix \ref{ap:lem}.\\
\qed

\subsection{Diffusion with small steps and long time}

In the case $x\in \R^N$,  the rescaling with small diffusion steps and long time, is given by
\beq
\label{SMeN}
\begin{cases}
\e \p_t n_\e (x,t) = \int_0^\infty\int_{\nu\in S^{N-1}}\left( n_\e(x+(e^{\e k}-1)\nu,t)  - n_\e(x,t) \right) \f{e^kdS dk}{|e^k-1|^{1+\al}} +n_\e(x,t) \,R(n_\e,I_\e)(x,t),\\
n_e(x,0)=n_\e^0(x),
\end{cases}
\eeq
with 
$$
I_\e(t)=\int n_\e(x,t)dx,
$$
Note that, in the case $N=1$, we retrieve \fer{SMe}. Replacing \fer{SMe} by \fer{SMeN}, and assumption \fer{as:biz} by
\beq
\label{as:bizN}
u_\e^0(x+h\nu) \leq u_\e^0(x)+A \log \left(1+|h|\right),\qquad \text{for all $x\in \R^N$, $h\in \R^+$ and $\nu\in S^{N-1}$},
\eeq
Theorems \ref{th:main} and \ref{th:main2} hold true for $x\in \R^N$: 

\begin{theorem}\label{th:mainN}
Let $x\in \R^N$ and $n_\e$ be the solution of \fer{SMeN} with \fer{re-1} and $u_\e=\e\log n_\e$. Assume \fer{as:maxR}, \fer{as:monR}, \fer{as:I}, \fer{as:ue0}, \fer{as:u0} and  \fer{as:bizN}. (i) Then, as $\e\to 0$,  $(I_\e)_\e$ converges locally uniformly to $I_0$ and $(u_\e)_\e$ converges locally uniformly to a  continuous function $u$ which is  Lipschitz continuous with respect to $x$ and continuous in $t$. Moreover, $u$ is the unique viscosity solution to the following equation
$$
\begin{cases}
 \p_t u-\int_0^\infty \int_{\nu\in S^{N-1}}\left( e^{ k D_x  u\cdot \nu }-1 \right) \f{e^k \,dS\,dk}{|e^k-1|^{1+\al}} =0,\\
u(x,0)=u^0(x),
\end{cases}
$$
and
$$
\|D_x u \|_{L^\infty(\R^N \times \R^+)}\leq A,\qquad \max_{x\in \R} u(x,t)=0.
$$
(ii) Finally, along subsequences as $\e\to 0$, $n_\e$ converges in $L^\infty\left( w*(0,\infty) ; \mathcal{M}^1(\R^N) \right)$ to a measure $n$, such that, $supp\, n \subset \{(x,t) \,|\, u(x,t)=0\}$. 
 \end{theorem}

\begin{theorem}\label{th:main2N}
Let $x\in \R^N$ and $n_\e$ be the solution of  \fer{SMe} with \fer{re-2}  and $u_\e=\e\log n_\e$. Assume \fer{as:ue0}, \fer{as:u0}, \fer{as:bizN} and \fer{as:n01}.

 (i) Then, as $\e\to 0$, $(u_\e)_\e$ converges locally uniformly to a  function $u$ that is Lipschitz continuous with respect to  $x$ and continuous in $t$. Moreover, $u$ is the viscosity solution to the following Hamilton-Jacobi equation
$$
\begin{cases}
\max \left( \p_t u-\int_0^\infty \int_{\nu\in S^{N-1}}\left( e^{ k D_x  u\cdot \nu }-1 \right) \f{e^k dSdk}{|e^k-1|^{1+\al}} -1 ,u\right)=0,\\
u(x,0)=u^0(x),
\end{cases}
$$
and
$$
 \|D_x u \|_{L^\infty(\R^N\times \R^+)}\leq A.
$$
 (ii) Moreover, as $\e\to 0$,
$$
\begin{cases}
n_\e\to 0, & \text{locally uniformly in  $\{(x,t)\in \R^N\times (0,\infty) \, |\, u(t,x)<0 \}$,}\\
 n_\e \to 1, & \text{locally uniformly in  $\mathrm{Int}\;\{(x,t)\in \R^N\times (0,\infty) \, |\,u(t,x)=0 \}$.}
\end{cases}
$$

\end{theorem}

\proof[Proof of Theorems \ref{th:mainN} and \ref{th:main2N}] The proofs of Theorems  \ref{th:main} and \ref{th:main2} can be easily adapted to prove 
Theorems  \ref{th:mainN} and \ref{th:main2N}. We only show the differences in the arguments for the regularity estimates. The remaining parts of the proofs are similar to the one-dimensional case. \\

\noindent 
(i) {\bf Uniform bounds from above and below}. Same type of inequalities as in \fer{boundue} and \fer{boundue2} can be proved for the equations above. By analogy to the proofs of  Theorems  \ref{th:reg} and \ref{th:reg2}, the key point  is to show that the following integral
$$
S=\displaystyle\int_{k\geq 0}\int_{\nu\in S^{N-1}} \left[ 
\dfrac { \left( |x|^2+1 \right)^{\f{A}{2\e} }}  {\left( |x+\nu(e^{\e k}-1)|^2+1 \right)^{\f{A}{2\e} } } -1\right] \f{e^k}{|e^k-1|^{1+\al}}\,dSdk
$$
is bounded. We show how this can be proved. The other parts of the proofs are similar.
\\

\noindent We split the integral term above  to two parts
$$
S=\int_0^\infty \int_{\nu\in S^{N-1},\,\nu\cdot e_1>0}  \left[ 
\dfrac { \left( |x|^2+1 \right)^{\f{A}{2\e} }}  {\left( |x+\nu(e^{\e k}-1)|^2+1 \right)^{\f{A}{2\e} } }+ \dfrac { \left( |x|^2+1 \right)^{\f{A}{2\e} }} {\left( |x-\nu(e^{\e k}-1)|^2+1 \right)^{\f{A}{2\e} }}  -2\right] \f{e^k}{|e^k-1|^{1+\al}}\,dSdk.$$
Note that 
$$
\left| |x| -(e^{\e k}-1) \right|  \leq | x+(e^{\e k}-1)\nu | \leq \left| |x| +(e^{\e k}-1) \right|,
$$
and
$$
\left| |x| -(e^{\e k}-1) \right|  \leq | x-(e^{\e k}-1)\nu | \leq \left| |x| +(e^{\e k}-1) \right|.
$$
Using the above inequalities and following the arguments in the proof of Theorem  \ref{th:reg} we obtain that, for  a large positive constant $C_N$, 
$$
\left| \displaystyle\int_{k\geq 1} \int_{\nu\in S^{N-1},\,\nu\cdot e_1>0}  \left[ 
\dfrac { \left( |x|^2+1 \right)^{\f{A}{2\e} }}  {\left( |x+\nu(e^{\e k}-1)|^2+1 \right)^{\f{A}{2\e} } }+ \dfrac { \left( |x|^2+1 \right)^{\f{A}{2\e} }} {\left( |x-\nu(e^{\e k}-1)|^2+1 \right)^{\f{A}{2\e} }}  -2\right] \f{e^k}{|e^k-1|^{1+\al}}\,dSdk
\right| \leq \f 1 2C_N.$$
To control the remaining part of the integral, that is 
$$
\left| \displaystyle\int_{0\leq k\leq 1} \int_{\nu\in S^{N-1},\,\nu\cdot e_1>0}  \left[ 
\dfrac { \left( |x|^2+1 \right)^{\f{A}{2\e} }}  {\left( |x+\nu(e^{\e k}-1)|^2+1 \right)^{\f{A}{2\e} } }+ \dfrac { \left( |x|^2+1 \right)^{\f{A}{2\e} }} {\left( |x-\nu(e^{\e k}-1)|^2+1 \right)^{\f{A}{2\e} }}  -2\right]  \f{e^k}{|e^k-1|^{1+\al}}\,dSdk
\right|,
$$
we first fix $\nu$, then use a Taylor expansion as in the proof of Theorem  \ref{th:reg}. Finally we integrate in $\nu$, to obtain, 
$$
\left| \displaystyle\int_{0\leq k\leq 1} \int_{\nu\in S^{N-1},\,\nu\cdot e_1>0}  \left[ 
\dfrac { \left( |x|^2+1 \right)^{\f{A}{2\e} }}  {\left( |x+\nu(e^{\e k}-1)|^2+1 \right)^{\f{A}{2\e} } }+ \dfrac { \left( |x|^2+1 \right)^{\f{A}{2\e} }} {\left( |x-\nu(e^{\e k}-1)|^2+1 \right)^{\f{A}{2\e} }}  -2\right]  \f{e^k}{|e^k-1|^{1+\al}}dSdk
\right| \leq \f 1 2 C_N.
$$
Combining the above arguments we obtain that $S$ is bounded.
\\

\noindent
(ii) \textbf{Logarithmic growth of $u_\e$.} 
We prove that
\beq
\label{log-esN}
u_\e(x+h\nu,t) \leq u_\e(x,t)+A \log \left(1+|h|\right),\qquad \text{for all $x\in \R^N$, $t\in\R^+$, $h\in \R^+$ and $\nu\in S^{N-1}$}.
\eeq

\noindent
For all $h\in \R$, $\nu \in S^{N-1}$ and $\e>0$, we define 
$$
w_{\e,h,\nu}(x,t)=u_\e(x+h \nu,t)-u_\e(x,t),\quad \text{for $t\geq 0$ and $x\in \R^N$.}
$$
We then compute 
$$
\begin{array}{rl}
\p_t w_{\e,h,\nu}(x,t) &=\displaystyle\int_{k\geq 0} \displaystyle\int_{\nu'\in S^{N-1}}\left[ e^{\f{u_\e(x+h \nu+(e^{\e k}-1)\nu',t)-u_\e(x+h\nu,t)}{\e}} - e^{\f{u_\e(x+(e^{\e k}-1)\nu',t)-u_\e(x,t)}{\e}}   \right]  \f{e^k dSdk}{|e^k-1|^{1+\al}}.
\end{array}
$$
Using a convexity inequality
 as before,
we deduce that
$$
\begin{array}{rl}
\p_t w_{\e,h,\nu}(x,t) &\leq\int_{k\geq 0} \int_{\nu'\in S^{N-1}} \left[ e^{\f{u_\e(x+h\nu+(e^{\e k}-1)\nu',t)-u_\e(x+h\nu,t)}{\e}} \left(\f{w_{\e,h,\nu}(x+(e^{\e k}-1)\nu',t)-w_{\e,h,\nu}(x,t)}{\e}\right)   \right]  \f{e^k dS\, dk}{|e^k-1|^{1+\al}}.
\end{array}
$$
Therefore, by the maximum principle and  \fer{as:biz} we obtain that for all $t>0$, $\e>0$ and $h,\,x\in\R$,
$$
w_{\e,h,\nu}(x,t)\leq \sup_x\, w_{\e,h,\nu}(x,0) \leq A\log (1+|h|),
$$
and hence \fer{log-esN} follows.\\

\appendix

\section{The proofs of Lemma \ref{lem:ex} and Lemma \ref{lem:exN}}
\label{ap:lem}
\subsection{The proof of Lemma \ref{lem:ex}}
In this section, we prove  Lemma \ref{lem:ex}. To this end, we let $\da<\f 1 2$ be a positive constant and suppose that $x>0$. The case with $x<0$ can be studied following similar arguments. We compute
$$
\begin{array}{rl}
\left| \f{(-\Delta )^{\f \al 2}g(x)}{g(x)} \right|&= 
\left| \displaystyle\int_0^\infty \left( \f{1+|x|^{1+\al}}{1+|x+h|^{1+\al}} +  \f{1+|x|^{1+\al}}{1+|x-h|^{1+\al}} -  2 \right) \f{dh}{|h|^{1+\al}} \right|\\
&\leq 
 \left| \displaystyle\int_{\R^+\setminus [0,\da]\cup [(1-\da)x,(1+\da)x]} \left( \f{1+|x|^{1+\al}}{1+|x+h|^{1+\al}} +  \f{1+|x|^{1+\al}}{1+|x-h|^{1+\al}} -  2 \right) \f{dh}{|h|^{1+\al}} \right|\\
&+\left| \displaystyle\int_{(1-\da)x\vee \da}^{(1+\da)x\vee \da} \left( \f{1+|x|^{1+\al}}{1+|x+h|^{1+\al}} +  \f{1+|x|^{1+\al}}{1+|x-h|^{1+\al}} -  2 \right) \f{dh}{|h|^{1+\al}} \right|\\
&+ \left| \displaystyle\int_0^\da \left( \f{1+|x|^{1+\al}}{1+|x+h|^{1+\al}} +  \f{1+|x|^{1+\al}}{1+|x-h|^{1+\al}} -  2 \right) \f{dh}{|h|^{1+\al}} \right|\   = \ I_1+I_2+I_3.
\end{array}
$$
We first notice that by easy computations one can obtain
$\
I_1 \leq \f{C}{\da^{(1+2\al)}}.
$
To control the second integral we write
$$
I_2 \leq \int_{(1-\da)x\vee \da}^{(1+\da)x\vee \da} \left( C+\f{1+|x|^{1+\al}}{1+|x-h|^{1+\al}} \right) \f{dh}{|h|^{1+\al}}=\int_{(1-\da)x\vee \da}^{(1+\da)x \vee \da} \left( C+\f{1}{|x|^{-(1+\al)}+\left(\f{|x-h|}{|x|}\right)^{1+\al}} \right) \f{dh}{|h|^{1+\al}}.
$$
Letting $\mu$ be an arbitrary small positive constant, we then use the Young's inequality to obtain that there exists a  positive constant $C$ such that
$$
\f{1}{|x|^{-(1+\al)}+\left( \f{|x-h|}{|x|} \right)^{1+\al} } \leq \f{C}{ |x|^{-(\mu+\al)} \left( \f{|x-h|}{|x|} \right)^{1-\mu} }=
C\f{|x|^{1+\al}}{|x-h|^{1-\mu}}.
$$
and hence,
$$\begin{array}{rl}
I_2 
&\leq C \displaystyle\int_{(1-\da)x\vee \da}^{(1+\da)x\vee \da} \left( 1+\f{|x|^{1+\al}}{|x-h|^{1-\mu}} \right) \f{dh}{|h|^{1+\al}}\ 
\leq \ C\f{1}{|\da|^\al}+C \displaystyle\int_{(1-\da)x}^{(1+\da)x} \f{1}{|x-h|^{1-\mu}} dh\leq
 C\left( \f{1}{|\da|^\al}+(\da |x|)^\mu \right).
 \end{array}
$$
Since this is true for arbitrarily small $\mu$ we obtain that 
$$
I_2 
\leq  C\left( \f{1}{|\da|^\al}+1 \right).
$$
To control $I_3$, we define
$$
f(x,h)=\f{1+|x|^{1+\al}}{1+|x+h|^{1+\al}}.
$$
We compute
$$ 
\f {\p }{\p h} f(x,h) = -(1+\al)\f{|x+h|^\al \left( 1+|x|^{1+\al} \right)}{ \left( 1+|x+h|^{1+\al} \right)^2}.
$$
It is easy to verify that for all $\eta_1,\, \eta_2\in [0,h]$, 
$$
|\f {\p }{\p h} f(x,\eta_1) - \f {\p }{\p h} f(x,-\eta_2) | \leq C |h|^{\al},
$$
for some constant $C$ independent of $|x|$ and $h$. It follows that
$$ 
| f(x+h) + f(x-h) -2 | \leq  C|h|^{1+\al}.
$$
and hence
$\ I_3 \leq C\da.
$

\noindent 
Fixing $0<\da < \f 1 2$, and combining the above inequalities we obtain that there exists a positive constant $C$ independent of $x$, such that
$$
\left| \f{(-\Delta )^{\f \al 2}g(x)}{g(x)} \right|\leq C.
$$

\subsection{The proof of Lemma \ref{lem:exN}}

Note that Lemma \ref{lem:exN} is the generalization of Lemma \ref{lem:ex} to the multidimensional case. We show that this generalization can be done easily. \\

\noindent
To this end, we split \fer{int-N}  to two parts
$$
\int_0^\infty \int_{\nu\in S^{N-1},\,\nu\cdot e_1>0} \left( \f{1+|x|^{1+\al}}{1+|x+h\nu|^{1+\al}} +\f{1+|x|^{1+\al}}{1+|x-h\nu|^{1+\al}}-2\right)\f{dS\,dh}{|h|^{1+\al}}.
$$
Note that 
$$
\left| |x| -h \right|  \leq | x+h\nu | \leq \left| |x| +h \right|,\qquad \left| |x| -h \right|  \leq | x-h\nu | \leq \left| |x| +h \right|.
$$
We fix  $0<\da < \f 1 2$ as in the proof of Lemma \ref{lem:ex}. Then, using the above inequalities and following the arguments in the proof of Lemma \ref{lem:ex} we obtain that, for  a large positive constant $C_N$ independent of $x$, 
$$
\left| \int_\da^\infty \int_{\nu\in S^{N-1},\,\nu\cdot e_1>0} \left( \f{1+|x|^{1+\al}}{1+|x+h\nu|^{1+\al}} +\f{1+|x|^{1+\al}}{1+|x-h\nu|^{1+\al}}-2\right)\f{dS\,dh}{|h|^{1+\al}} \right|
\leq \f 1 2 \,C_N.
$$

\noindent
To control the remaining term of the integral, that is 
$$
\int_0^\da \int_{\nu\in S^{N-1},\,\nu\cdot e_1>0} \left( \f{1+|x|^{1+\al}}{1+|x+h\nu|^{1+\al}} +\f{1+|x|^{1+\al}}{1+|x-h\nu|^{1+\al}}-2\right)\f{dS\,dh}{|h|^{1+\al}},
$$
 we first fix $\nu$, then do the same computation as in the proof of Lemma \ref{lem:ex}. Finally we integrate in $\nu$, to obtain, 
$$
\int_0^\da \int_{\nu\in S^{N-1},\,\nu\cdot e_1>0} \left( \f{1+|x|^{1+\al}}{1+|x+h\nu|^{1+\al}} +\f{1+|x|^{1+\al}}{1+|x-h\nu|^{1+\al}}-2\right)\f{dS\,dh}{|h|^{1+\al}}\leq \f 1 2 \,C_N.
$$
Combining the above arguments we obtain \fer{gN}.

\bigskip

%
%
%
%

\textbf{Acknowledgments:} 
The second Author is very grateful to Guy Barles and Jean-Michel Roquejoffre for  fruitful discussions. She wishes  also to acknowledge partial support by the french ANR projects KIBORD ANR-13-BS01-0004 and MODEVOL ANR-13-JS01-0009. The two authors acknowledge partial support  by the  Chaire Mod\'elisation Math\'ematique et Biodiversit\'e VEOLIA-\'Ecole Polytechnique-MNHN-F.X. and  the  ANR project MANEGE
ANR-09-BLAN-0215.
%

\begin{thebibliography}{10}

\bibitem{Baeumer}
B.~Baeumer, M.~Kovacs, and M.M.~Meerschaert.
\newblock Fractional reproduction-dispersal equations and heavy tail dispersal kernels.
\newblock {\em Bull. Math. Biol.}, 69:2281--2297, 2007.

\bibitem{GB:94}
G.~Barles.
\newblock {\em Solutions de viscosit\'e des \'equations de
  {H}amilton-{J}acobi}, volume~17 of {\em Math\'ematiques \& Applications
  (Berlin) [Mathematics \& Applications]}.
\newblock Springer-Verlag, Paris, 1994.

\bibitem{GB.LE.PS:90}
G.~Barles, L.~C. Evans, and P.~E. Souganidis.
\newblock Wavefront propagation for reaction-diffusion systems of {PDE}.
\newblock {\em Duke Math. J.}, 61(3):835--858, 1990.



\bibitem{GB.CI:08}
G.~Barles and C.~Imbert.
\newblock Second-order elliptic integro-differential equations: Viscosity
  solutionsÕ theory revisited.
\newblock {\em Ann. Inst. H. Poincar\'e Anal. Non Lin\'eaire}, 25:567--585,
  2008.


\bibitem{GB.SM.BP:09}
G.~Barles, S.~Mirrahimi, and B.~Perthame.
\newblock Concentration in {L}otka-{V}olterra parabolic or integral equations:
  a general convergence result.
\newblock {\em Methods Appl. Anal.}, 16(3):321--340, 2009.

\bibitem{GB.BP:88}
G.~Barles and B.~Perthame.
\newblock Exit time problems in optimal control and vanishing viscosity method.
\newblock {\em SIAM J. Control Optim.}, 26(5):1133--1148, 1988.

\bibitem{HB.JR.LR:11}
H.~Berestycki, J.-M. Roquejoffre, and L.~Rossi.
\newblock The periodic patch model for population dynamics with fractional
  diffusion.
\newblock {\em Discrete Contin. Dyn. Syst. Ser. S}, pages 1--13, 2011.


\bibitem{CB.EC:09}
C.~Br\"{a}ndle and E.~Chasseigne.
\newblock Large deviations estimates for some non-local equations {I}.~fast
  decaying kernels and explicit bounds.
\newblock {\em Nonlinear Anal.}, (11):5572--5586, 2009.

\bibitem{CB.EC:13}
C.~Br\"{a}ndle and E.~Chasseigne.
\newblock Large deviations estimates for some non-local equations. general
  bounds and applications.
\newblock {\em Trans. Amer. Math. Soc.}, pages 3437--3476, 2013.




\bibitem{XC.AC.JR:12}
X.~Cabr\'e, A.-C. Coulon, and J.-M. Roquejoffre.
\newblock Propagation in {Fisher-KPP} type equations with fractional diffusion
  in periodic media.
\newblock {\em C. R. Math. Acad. Sci. Paris}, 350:885--890, 2012.

\bibitem{JR.XC:09}
X.~Cabr\'e and J.-M. Roquejoffre.
\newblock Propagation de fronts dans les Žquations de {F}isher-{KPP} avec
  diffusion fractionnaire.
\newblock {\em C.R. Acad. Sci. Paris}, 347:1361--1366, 2009.


\bibitem{XC.JR:13}
X.~Cabr\'e and J.-M. Roquejoffre.
\newblock The influence of fractional diffusion on {F}isher-{KPP} equations.
\newblock {\em Comm. Math. Phys}, 320:679--722, 2013.



\bibitem{AC.JR:12}
A.C. Coulon and J.M. Roquejoffre.
\newblock Transition between linear and exponential propagation in
  {F}isher-{KPP} type reaction-diffusion equations.
\newblock {\em Comm. Partial Differential Equations}, 37:2029--2049, 2012.


\bibitem{C.I.L:92}
M.~G. Crandall, H.~Ishii, and P.-L. Lions.
\newblock User's guide to viscosity solutions of second order partial
  differential equations.
\newblock {\em Bull. Amer. Math. Soc. (N.S.)}, 27(1):1--67, 1992.

\bibitem{OD.PJ.SM.BP:05}
O.~Diekmann, P.-E. Jabin, S.~Mischler, and B.~Perthame.
\newblock The dynamics of adaptation: an illuminating example and a
  {H}amilton-{J}acobi approach.
\newblock {\em Th. Pop. Biol.}, 67(4):257--271, 2005.




\bibitem{HE:10}
H.~Engler.
\newblock On the speed of spread for fractional reaction-diffusion equations.
\newblock {\em Int. J. Differ. Equ.},  Article ID 315421, 16 p., 2010.


\bibitem{LE.PS:89}
L.~C. Evans and P.~E. Souganidis.
\newblock A {PDE} approach to geometric optics for certain semilinear parabolic
  equations.
\newblock {\em Indiana Univ. Math. J.}, 38(1):141--172, 1989.

\bibitem{MFb:85}
M~Freidlin.
\newblock {\em Functional integration and partial differential equations},
  volume 109 of {\em Annals of Mathematics Studies}.
\newblock Princeton University Press, Princeton, NJ, 1985.

\bibitem{MF:85}
M~Freidlin.
\newblock Limit theorems for large deviations and reaction-diffusion equations.
\newblock {\em The Annals of Probability}, 13(3):639--675,  1985.



\bibitem{Gurney}
W.S.~Gurney and R.M.~Nisbet.
\newblock The regulation of inhomogeneous populations.
\newblock {\em J. Theor. Biol.}, 52:441--457,1975.

\bibitem{Imkeller}
P.~Imkeller, I.~Pavlyukevich.
\newblock First exit times of SDEs driven by stable L\'evy processes.  
\newblock {\em Stoch. Process. Appl. }, 116 (4):611--642, 2006.

\bibitem{Jourdain}
B.~Jourdain, S.~M\'el\'eard, and W.~Woyckynski.
\newblock L\'evy flights in evolutionary ecology.
\newblock {\em J. Math. Biol.}, 65:677--707, 2012.


\bibitem{JG:11}
J.~Garnier.
\newblock Accelerating solutions in integro-differential equations.
\newblock {\em SIAM Journal on Mathematical Analysis}, (4):1955--1974, 2011.


\bibitem{GB.BP:08}
B.~Perthame and G.~Barles.
\newblock Dirac concentrations in {L}otka-{V}olterra parabolic {PDE}s.
\newblock {\em Indiana Univ. Math. J.}, 57(7):3275--3301, 2008.


\bibitem{Sato}
K.~Sato.
\newblock \emph{L\'evy processes and Infinitely Divisible Distributions}.
\newblock { Cambridge Studies in Advanced Math.}, 68, Cambridge University Press, 1999.

\end{thebibliography}

%

\end{document}